\documentclass[12pt]{article}
\usepackage[top=72pt,bottom=96pt,left=72pt,right=72pt]{geometry}
\usepackage{mathrsfs}
\usepackage{amssymb}
\def\1{\mathbf{1}}

\def\A{\mathscr{A}}
\def\aux{\mathrm{aux}}
\def\bm{\circledast}
\def\C{\mathbb{C}}
\def\End{\mathrm{End}\,}
\def\Ext{\mathrm{Ext}\,}
\def\Fr{\mathrm{Frame}}
\def\g{\mathfrak{g}}
\def\gl{\mathfrak{gl}}
\def\GL{\mathbf{GL}}
\def\h{\mathfrak{h}}

\def\Hom{\mathrm{Hom}\,}
\def\id{\mathrm{id}}
\def\im{\mathbf{im}}
\def\Jac{\mathbf{Jac}}
\def\Jet{\mathrm{Jet}}
\def\k{\mathfrak{k}}
\def\ker{\mathbf{ker}}
\def\L{\Lambda}
\def\max{\mathrm{max}}
\def\mod{\mathrm{mod}}
\def\M{\mathfrak{M}}
\def\n{\mathfrak{n}}
\def\N{\mathbb{N}}
\def\pr{\mathrm{pr}\,}
\def\qed{\hfill\ensuremath{\Box}}
\def\R{\mathbb{R}}
\def\rep{\mathrm{rep}}
\def\S{\mathrm{Sym}\,}
\def\Singer{\mathrm{Singer}}
\def\so{\mathfrak{so}}
\def\span{\mathrm{span}}
\def\SU{\mathbf{SU}}
\def\stab{\mathfrak{stab}\,}

\def\tf{\mathrm{tf}}
\def\tfrac{\textstyle\frac}
\def\tr{\mathrm{tr}}
\def\X{\mathfrak{X}}
\def\Y{\mathfrak{Y}}
\def\Z{\mathbb{Z}}
\newtheorem{Proposition}{Proposition}[section]
\newtheorem{Lemma}[Proposition]{Lemma}
\newtheorem{Theorem}[Proposition]{Theorem}
\newtheorem{Corollary}[Proposition]{Corollary}
\newtheorem{Definition}[Proposition]{Definition}

\newtheorem{Remark}[Proposition]{Remark}
\newenvironment{Reference}[1]{\pfill\textbf{#1} \textit\bgroup}{\egroup\par}
\def\proof{\noindent\textbf{Proof:}\quad}
\def\pfill{\par\vskip6pt plus2pt minus3pt\noindent}
\begin{document}
\def\today{July 20th, 2017}
\title{Moduli Spaces of Affine Homogeneous Spaces}
\author{Gregor Weingart\footnote{Instituto de Matem\'aticas (Cuernavaca),
 Universidad Nacional Aut\'onoma de M\'exico, Avenida Universidad s/n,
 Colonia Lomas de Chamilpa, 62210 Cuernavaca, Morelos, MEXICO;
 \texttt{gw@matcuer.unam.mx}}}
\maketitle
\begin{center}
 \textbf{Abstract}
 \\[11pt]
 \parbox{412pt}{%
  Apart from global topological problems an affine homogeneous space
  is locally described by its curvature, its torsion and a slightly
  less tangible object called its connection in a given base point. Using
  this description of the local geometry of an affine homogeneous space
  we construct an algebraic variety $\M(\,\gl\,V\,)$, which serves as a
  coarse moduli space for the local isometry classes of affine homogeneous
  spaces of dimension $\dim\,V$. Moreover we associate a $\S\,V^*$--comodule
  to a point in $\M(\,\gl\,V\,)$ and use its Spencer cohomology in order
  to describes the infinitesimal deformations of this point in the true
  moduli space $\M_\infty(\,\gl\,V\,)$.}
 \\[11pt]
 \textbf{MSC2010:\quad 53C30;\ 22F30}
\end{center}
\section{Introduction}
\label{intro}
 Homogeneous spaces comprise a class of smooth manifolds of particular
 interest in differential geometry, because many geometric calculations
 reduce essencially to linear algebra in the presence of a transitively
 acting Lie group. In particular all of the geometry of a homogeneous
 space except for its global topology can be represented in terms of
 converging formal power series on a formal neighborhood of the base
 point. This representation of the local geometry of a homogeneous space
 by convergent power series will be used implicitly in this article to
 construct moduli spaces of isometry classes of and a corresponding
 deformation theory for affine homogeneous spaces with or without
 additional geometric structures.

 A particular advantage of our approach in comparison to the moduli spaces
 of locally homogeneous spaces constructed in \cite{tsu} is that the covariant
 derivatives of the curvature and/or the torsion in the construction of the
 moduli spaces are replaced by the Christoffel symbols of the connection
 in the base point. The locally affine homogeneous spaces considered in
 \cite{sing}\cite{kov}\cite{bo}\cite{tric}\cite{tsu} for example can be thought
 of as affine homogeneous spaces with variable Christoffel symbols or
 equivalently as smooth maps to the moduli spaces constructed below, the
 derivatives of these maps can thus be studied using the deformation
 theory of affine homogeneous spaces arising naturally in our approach.
 Put differently the examples of locally, but not globally affine homogeneous
 spaces owe their existence to the non--vanishing of certain Spencer
 cohomology spaces, which replace the Chevalley--Eilenberg cohomology
 spaces in the deformation theory of Lie algebras.

 \pfill
 Representing the local geometry of an affine homogeneous spaces in terms
 of converging formal power series obliterates its global topology as a
 manifold, hence we will essentially consider only formal affine homogeneous
 spaces defined as a pair $\g\,\supset\,\h$ of Lie algebras endowed with an
 $\h$--equivariant linear map $A:\,\g\longrightarrow\End\,\g/\h$ extending
 the isotropy representation $\star:\,\h\longrightarrow\End\,\g/\h$. This
 algebraic simplicification is bought at a prize however: The integration of
 a formal affine homogeneous space into a real manifold will fail in general
 to produce an affine homogeneous space, because the Lie algebra $\h$ may
 not integrate to a closed subgroup of the simply connected Lie group $G$
 corresponding to the Lie algebra $\g$. Nevertheless every formal affine
 homogeneous space corresponds to a homogeneous transitive Lie algebroid
 over a homogeneous space endowed with a left invariant connection as
 discussed briefly in the comments after Definition \ref{fahs}.

 In order to construct the moduli spaces of formal affine homogeneous spaces
 we fix a model vector space $V$ of the appropiate dimension and augment a
 given formal affine homogeneous space $\g\,\supset\,\h$ endowed with
 $A:\,\g\longrightarrow\End\,\g/\h$ by an isomorphism or frame $F:\,
 V\longrightarrow\g/\h$ and a not necessarily $\h$--equivariant split
 $\g/\h\longrightarrow\g$ of the canonical projection. With these two
 additional pieces of data in place we may associate a
 connection--curvature--torsion triple
 $$
  (\;A,\;R,\;T\;)\;\;\in\;\;
  V^*\,\otimes\,\End\,V
  \;\times\;\L^2V^*\,\otimes\,\End\,V
  \;\times\;\L^2V^*\,\otimes\,V
 $$
 to an augmented formal affine homogeneous space, which encodes the geometry
 of $\g/\h$ completely. Our main result characterizes the subset $\M(\,\gl\,
 V\,)$ of connection--curvature--torsion triples arising from augmented formal
 affine homogeneous spaces by means of a formally infinite, but actually
 finite system of explicit homogeneous algebraic equations:

 \begin{Reference}
  {Theorem \ref{mfh} (Algebraic Variety of Affine Homogeneous Spaces)}
  \hfill\break
  A connection--curvature--torsion triple $(A,\,R,\,T)$ on a vector space
  $V$ represents a formal affine homogeneous space $\g\,\supset\,\h$
  augmented by a frame isomorphism $F:\,V\longrightarrow\g/\h$ and a split
  $\g/\h\longrightarrow\g$ of the canonical projection, if and only if
  $(A,\,R,\,T)$ satisfies the formal first and second Bianchi identities
  $d^{(A,T)}T\,=\,R\wedge\id$ and $d^{(A,T)}R\,=\,0$ of degrees $2,\,3$
  respectively and for all $r\,\geq\,0$ the following homogeneous
  equations of degrees $r+3$ and $r+4$
  \begin{eqnarray*}
   \Big(\;Q(A,T)\;-\;R\;\Big)\,\bm\,\Big(\;\underbrace{A\,\bm\,(\,A\,\bm\,
   (\,\ldots\,(\,A}_{r\;\,\mathrm{times}}\,\bm\,T\,)\,\ldots\,)\,)\;\Big)
   &=&0
   \\
   \Big(\;Q(A,T)\;-\;R\;\Big)\,\bm\,\Big(\;\underbrace{A\,\bm\,(\,A\,\bm\,
   (\,\ldots\,(\,A}_{r\;\,\mathrm{times}}\,\bm\,R\,)\,\ldots\,)\,)\;\Big)
   &=&0
  \end{eqnarray*}
  where $Q(A,T)_{x,y}\,:=\,[\,A_x,\,A_y\,]\,-\,A_{A_xy-A_yx-T(x,y)}$ is
  the approximate curvature of $(A,\,R,\,T)$.
 \end{Reference}

 \pfill
 This result generalizes the well--known description of local isometry
 classes of symmetric spaces in terms of their curvature $R$ \cite{helg}
 and the classification of manifolds with parallel curvature and torsion
 by Ambrose and Singer \cite{as}. Eventually we will introduce the
 equivalence relation of infinite order contact $\sim_\infty$ on
 connection--curvature--torsion triples, which describes precisely the
 effect of changing the frame $F$ and/or the split $\g/\h\longrightarrow\g$.
 The true moduli space of isometry classes of formal affine homogeneous
 spaces is thus the quotient:
 $$
  \M_\infty(\;\gl\,V\;)\;\;:=\;\;\M(\;\gl\,V\;)/_{\displaystyle\sim_\infty}
 $$
 A subtle invariant of a connection--curvature--torsion triple $(A,\,R,\,T)
 \,\in\,\M(\,\gl\,V\,)$ plays a prominent role in the proof of Theorem
 \ref{mfh}, the stabilizer filtration of $\End\,V$ by subalgebras
 $$
  \End\,V\;=\;\ldots\;=\;\h_{-2}\;=\;\h_{-1}\;\supsetneq\;\h_0
  \;\supsetneq\;\ldots\;\supsetneq\;\h_{s-1}
  \;\supsetneq\;\h_s\;=\;\h_{s+1}\;=\;\ldots\;=\;\h_\infty
 $$
 introduced by Singer \cite{sing}, in particular the minimal $s\,\geq\,-1$
 with equality $\h_s\,=\,\h_{s+1}$ is nowadays called the Singer invariant
 of $(A,\,R,\,T)$. Our recursive Definition \ref{sfilt} of the stabilizer
 filtration differs significantly from Singer's definition akin to Lemma
 \ref{if} and is certainly easier to use in actual calculations. The graded
 vector space of sucessive filtration quotients in the stabilizer
 filtration associated to a connection--curvature--torsion triple
 $$
  \h^\bullet
  \;\;:=\;\;
  \bigoplus_{r\,\in\,\Z}\Big(\;\h_{r-1}/_{\displaystyle\h_r}\;\Big)
  \;\;=\;\;
  (\,\End\,V/_{\displaystyle\h_0}\,)\;\oplus\;\ldots\;\oplus\;
  (\,\h_{s-1}/_{\displaystyle\h_\infty}\,)
 $$
 is naturally a $\S\,V^*$--comodule and thus allows us to associate a
 cohomology theory to every point $(A,\,R,\,T)$ in the moduli space
 $\M(\,\gl\,V\,)$, namely the Spencer cohomology $H^{\bullet,\circ}(\,\h\,)$
 of the comodule $\h^\bullet$. The second main result of this article links
 the special Spencer cohomology spaces $H^{\bullet,1}(\,\h\,)$ to a geometric
 filtration on the formal tangent space to the true moduli space $\M_\infty
 (\,\gl\,V\,)$ in the point represented by the connection--curvature--torsion
 triple $(A,\,R,\,T)$:
 \begin{equation}\label{ftseq}
  H^{\bullet,\,1}(\;\h\;)
  \;\;=\;\;
  T_{[A,R,T]}\M^\bullet_\infty(A,R,T)
  /_{\displaystyle T_{[A,R,T]}\M^{\bullet-1}_\infty(A,R,T)}
 \end{equation}
 Philosophically this equality reflects the simple fact that a vector field
 representing an $r$--jet solution to the affine Killing equation, which can
 not be lifted to an $r+1$--jet solution, can {\em not} be an affine Killing
 field, hence its flow will change the underlying geometry, but gently enough
 to stay in contact with the original geometry up to order $r-2$. In a rather
 precise sense equation (\ref{ftseq}) is a quantitative version of Singer's
 Theorem \cite{sing}, which characterizes the globally among the locally
 Riemannian homogeneous spaces, because the Spencer cohomology spaces are
 trivial by construction for degrees $\bullet\,>\,s$ larger than the Singer
 invariant so that no further deformations are possible. A paradoxical aspect
 of the preceeding construction is that Spencer cohomology is usually
 introduced to describe the set of $r+1$--jet solutions as an affine space
 bundle over the set of $r$--jet solutions, however it works just as well
 for the affine Killing equation, where the $r$--jet solutions are an affine
 space over the $r+1$--jet solutions.

 \pfill
 In Section \ref{ginv} of this article we will recall many well--known
 facts about homogeneous spaces, in particular we will discuss a very
 convenient calculus for left invariant connections on homogeneous vector
 bundles. Section \ref{afs} describes two equivalent ways to encode the
 local geometry of a formal affine homogeneous space in an algebra endowed
 with a skew bracket, one of these algebras has been studied in \cite{nt}.
 In Section \ref{class} we introduce connection--curvature--torsion triples
 and derive the algebraic equations defining of the coarse moduli space
 $\M(\,\gl\,V\,)$ as an algebraic variety. Additional parallel geometric
 structures like Riemannian metrics or almost complex structures are added
 to the picture in Section \ref{par}, whereas the final Section \ref{fts}
 constructs the Spencer cohomology of a formal affine homogeneous space
 and calculates this cohomology for the family of examples of Riemannian
 homogeneous spaces with large Singer invariant constructed by C.~Meusers
 \cite{m}, which turns out to be a maximal family in the sense that it is
 closed under all its possible deformations.

 \pfill
 The framework of this article was developped in intensive collaboration
 with C.~Meusers, the recursive definition of the stabilizer filtration for
 example arose directly from these discussion and was subsequently published
 in \cite{m}. Besides C.~Meusers the author would like to thank W.~Ballmann
 for his support and encouragement in writing up this article.
\section{Left Invariant Connections}
\label{ginv}
 Certainly the most important single concept in studying the geometry
 of a homogeneous space $G/H$ is the notion of a homogeneous vector
 bundle over $G/H$, a vector bundle endowed with a left $G$--action
 on its total space, which covers the canonical left action of $G$
 on $G/H$ and is linear in every fiber. The additional $G$--action
 on the total space makes the vector space of sections of a homogeneous
 vector bundle a representation of the group $G$ generalizing the left
 regular representation of $G$ on $C^\infty(G)$. A differential operator
 between sections of homogeneous vector bundles commuting with the respective
 representations of $G$ is called for this reason a left invariant differential
 operator. In this section we will focus on the algebraic properties of a
 specific subclass of left invariant differential operators, left invariant
 connections. The detailed algebraic formalism introduced in this section
 to describe left invariant connections will turn out to be quite useful
 for our subsequent calculations.

 \pfill
 Recall that a homogeneous vector bundle on a homogeneous space $G/H$
 is a vector bundle on which $G$ acts from the left by vector bundle
 morphisms covering the left action of $G$ on $G/H$. General nonsense
 implies that the category of homogeneous vector bundles over $G/H$
 under $G$--equivariant vector bundle homomorphisms is equivalent to
 the category of representations of $H$ under $H$--equivariant linear
 maps. One possible choice for the $H$--representation in this
 correspondence is simply the fiber $\Sigma\,:=\,\Sigma_{eH}(G/H)$
 of the homogeneous vector bundle $\Sigma(G/H)$ over the base point
 $eH\,\in\,G/H$, which is a representation of $H$ by the very definition
 of a homogeneous vector bundle. For the tangent and cotangent bundles
 however these fibers are usually replaced by the representations $\g/\h$
 and $(\g/\h)^*$, which are isomorphic to the fibers $T_{eH}(G/H)$ and
 $T^*_{eH}(G/H)$ respectively via the $H$--equivariant isomorphism:
 $$
  \g/\h\;\stackrel\cong\longrightarrow\;T_{eH}(G/H),\qquad
  X\,+\,\h\;\longmapsto\;\left.\frac d{dt}\right|_0e^{tX}\,H
 $$

 \begin{Definition}[Homogeneous Vector Bundles]
 \hfill\label{hvb}\break
  A homogeneous vector bundle on a homogeneous space $G/H$ is a vector
  bundle $\Sigma(G/H)$ on $G/H$ endowed with a left $G$--action on its
  total space, which covers the canonical left action $\star:\,G\times G/H
  \longrightarrow G/H,\,(\gamma,gH)\longmapsto\gamma gH,$ of $G$ on $G/H$
  and is linear in every fiber:
  $$
   \star:\quad G\;\times\;\Sigma(\,G/H\,)
   \;\longrightarrow\;\Sigma(\,G/H\,),\qquad(\,\gamma,\,s\,)
   \;\longmapsto\;\gamma\,\star\,s
  $$
  The vector space of sections of a homogeneous vector bundle becomes a
  representation of $G$
  $$
   L:\quad G\;\times\;\Gamma\,\Sigma(\,G/H\,)
   \;\longrightarrow\;\Gamma\,\Sigma(\,G/H\,),
   \qquad(\;\gamma,\,s\;)\;\longmapsto\;L_\gamma s
  $$
  by means of $(L_\gamma s)(\,gH\,)\,:=\,\gamma\,\star\,s(\,\gamma^{-1}gH\,)$
  for all $\gamma\,\in\,G$ and all sections $s\,\in\,\Gamma\,\Sigma(G/H)$.
 \end{Definition}

 \pfill
 In the same vein a homogeneous fiber bundle over a homogeneous space $G/H$
 is a fiber bundle over $G/H$ endowed with a left action of $G$ on its total
 space, which covers the canonical left action on $G/H$. In particular the
 group $G$ itself can be considered as a homogeneous principal $H$--bundle
 over $G/H$, where the commuting left and right multiplications in
 $$
  G\,\times\,G\,\times\,H\;\longrightarrow\;G,\qquad
  (\,\gamma,\,g,\,h\,)\;\longmapsto\;\gamma\,g\,h
 $$
 define the left $G$--action on the total space and the principal $H$--bundle
 structure respectively. In turn a left invariant principal connection is a
 differential form $\omega\,\in\,\Gamma(\,T^*G\otimes\h\,)$ on $G$ with
 values in $\h$ invariant under the left $G$--action, such that the axiom
 for principal $H$--connections
 \begin{equation}\label{pca}
  \omega_{g_0\,h_0}(\;\left.\frac d{dt}\right|_0g_t\,h_t\;)
  \;\;=\;\;
  \mathrm{Ad}_{h_0^{-1}}\,\omega_{g_0}(\;\left.\frac d{dt}\right|_0g_t\;)
  \;+\;\left.\frac d{dt}\right|_0h_0^{-1}\,h_t
 \end{equation}
 is satisfied for all curves $t\longmapsto g_t$ in $G$ and $t\longmapsto h_t$
 in $H$. The Maurer--Cartan form $\theta\,\in\,\Gamma(\,T^*G\otimes\g\,)$
 provides a left invariant trivialization of $TG$, hence every left invariant
 connection is necessarily of the form $\omega\circ\theta$ for a linear map
 $\omega:\,\g\longrightarrow\h$. However the principal connection axiom
 (\ref{pca}) requires $\omega$ to be an $H$--equivariant section of the
 inclusion:

 \begin{Definition}[Left Invariant Principal Connections]
 \hfill\break
  A left invariant principal connection on $G$ considered as a homogeneous
  principal $H$--bundle over $G/H$ is an $H$--equivariant section
  $\omega:\;\g\,\longrightarrow\,\h$ of the short exact sequence:
  $$
   0\;\longrightarrow\;\h\;\longrightarrow\;\g\;
   \longrightarrow\;\g/\h\;\longrightarrow\;0
  $$
 \end{Definition}

 \pfill
 In consequence of this definition the set of left invariant principal
 connections on a given homogeneous space $G/H$ is either the empty set
 or an affine space modelled on $\Hom_H(\g/\h,\h)$. In fact the set of
 all left invariant principal connections on $G/H$ is precisely the
 preimage of $\id_\h\,\in\,\Hom_H(\h,\h)$ under the map $\Hom_H(\g,\h)
 \longrightarrow\Hom_H(\h,\h)$ in the long exact sequence
 $$
  0
  \;\longrightarrow\;
  \Hom_H(\,\g/\h,\,\h\,)
  \;\longrightarrow\;
  \Hom_H(\,\g,\,\h\,)
  \;\longrightarrow\;
  \Hom_H(\,\h,\,\h\,)
  \;\stackrel\delta\longrightarrow\;
  \Ext^1_H(\,\g/\h,\,\h\,)
  \;\longrightarrow\;
 $$
 associated to the short exact sequence of $H$--representations:
 $$
  0\;\longrightarrow\;\h\;
  \longrightarrow\;\g\;\longrightarrow\;\g/\h\;
  \longrightarrow\;0
 $$
 The extension class $\delta(\,\id_\h\,)\in\Ext^1_H(\,\g/\h,\,\h\,)$ is
 thus the unique obstruction against the exis\-tence of a left invariant
 principal connection on a homogeneous space $G/H$. A homogenous space
 $G/H$ with vanishing $\delta(\,\id_\h\,)\,=\,0$ is properly called a
 reductive homogeneous space, in the literature however this notion
 implicitly includes the choice of some left invariant principal connection
 $\omega:\,\g\longrightarrow\h$ or other. Concealing this {\em arbitrary}
 choice in the inconspi\-cuous adjective {\em reductive} it is somewhat
 duplicitous though to call the induced connections on homogeneous vector
 bundles {\em canonical} connections.

 With the cotangent bundle $T^*(G/H)$ of a homogeneous space $G/H$ being
 homogeneous the group $G$ acts not only on the vector space of sections
 $\Gamma(\,\Sigma(G/H)\,)$ of a homogeneous vector bundle $\Sigma(G/H)$,
 but also on $\Gamma(\,T^*(G/H)\otimes\Sigma(G/H)\,)$. In turn a linear
 connection
 \begin{equation}\label{cleft}
  \nabla:\quad\Gamma(\;\Sigma(G/H)\;)\;
  \longrightarrow\;\Gamma(\;T^*(G/H)\,\otimes\,\Sigma(G/H)\;)
 \end{equation}
 on $\Sigma(G/H)$ is called a left invariant connection on $\Sigma(G/H)$
 provided $\nabla$ is $G$--equivariant in the sense $L_\gamma(\nabla s)
 \,=\,\nabla(L_\gamma s)$ for all sections $s\,\in\,\Gamma(\,\Sigma
 (G/H)\,)$ and every $\gamma\,\in\,G$. For the purpose of this article
 we want to replace this definition by a more algebraic formulation:

 \begin{Definition}[Left Invariant Connections]
 \hfill\label{lconn}\break
  A left invariant connection on a homogeneous vector bundle $\Sigma(G/H)$ is
  an $H$--equivariant extension $A:\,\g\,\longrightarrow\,\End\,\Sigma$ of the
  infinitesimal representation $\star:\,\h\,\longrightarrow\,\End\,\Sigma$ of
  $\h$ on $\Sigma$.
 \end{Definition}

 \noindent
 Similar to the classification of left invariant principal connections
 the set of left invariant connections $A:\,\g\longrightarrow\End\,\Sigma$
 on $\Sigma(G/H)$ is the preimage of the infinitesimal representation
 $\star:\,\h\longrightarrow\End\,\Sigma$ under $\Hom_H(\g,\End\,\Sigma)
 \longrightarrow\Hom_H(\h,\End\,\Sigma)$ in the long exact sequence:
 $$
  \ldots\;\longrightarrow\;
  \Hom_H(\,\g,\,\End\,\Sigma\,)
  \;\longrightarrow\;
  \Hom_H(\,\h,\,\End\,\Sigma\,)
  \;\stackrel\delta\longrightarrow\;
  \Ext^1_H(\,\g/\h,\,\End\,\Sigma\,)
  \;\longrightarrow\;\ldots
 $$
 In consequence the existence of left invariant connections $A:\,\g
 \longrightarrow\End\,\Sigma$ on $\Sigma(G/H)$ is obstructed by the
 extension class $\delta(\,\star\,)\,\in\,\Ext^1_H(\g/\h,\End\,\Sigma)$,
 its vanishing provides us with an affine space of left invariant connections
 on $\Sigma(G/H)$ modelled on $\Hom_H(\g/\h,\End\,\Sigma)$. Evidently
 every left invariant principal connection $\omega:\,\g
 \longrightarrow\h$ induces a left invariant connection $A^\omega:
 \,\g\longrightarrow\End\,\Sigma$ on every homogeneous vector bundle
 $\Sigma(G/H)$ by means of $A^\omega(X)\,:=\,\omega(X)\,\star$. Due to
 this universality every non--vanishing obstruction $\delta(\,\star\,)
 \,\neq\,0$ for some $H$--representation $\Sigma$ is of course an
 obstruction against the existence of a left invariant principal
 connection as well.

 For the moment we will not discuss the precise relationship between
 the analytic and the algebraic definition of a left invariant connection
 on a homogeneous vector bundle $\Sigma(G/H)$ over a homogeneous space
 $G/H$. Instead we want to point out a couple of very remarkable analogies,
 which taken together should leave no doubt that these two concepts are
 intimately related. The algebraic Definition \ref{lconn} for example is
 easy to use with constructions of linear algebra like dual spaces, direct
 sums and tensor products, say the dual left invariant connection on
 $\Sigma^*$ is simply given by $-A^*:\,\g\longrightarrow\End\,\Sigma^*,\,
 X\longmapsto-A_X^*,$ whereas $A^\Sigma\oplus A^{\tilde\Sigma}$ and
 $A^\Sigma\otimes\id\,+\,\id\otimes A^{\tilde\Sigma}$ are left invariant
 connections on the direct sum $\Sigma\oplus\tilde\Sigma$ and tensor product
 $\Sigma\otimes\tilde\Sigma$ of two representations $\Sigma,\,\tilde\Sigma$
 endowed with left invariant connections $A^\Sigma,\,A^{\tilde\Sigma}$.

 \begin{Lemma}[Left Invariant Connections and Rigidity]
 \hfill\label{lir}\break
  If there exists a left invariant connection on the tangent bundle of
  a homogeneous space $G/H$, then the kernel $\n\,:=\,\ker\,\star$ of
  the isotropy representation $\star:\,\h\longrightarrow\End\,\g/\h,
  \,Z\longmapsto Z\,\star,$ is an ideal not only in the isotropy algebra
  $\h$, but already in the full Lie algebra $\g$.
 \end{Lemma}

 \proof
 According to the preceeding discussion a left invariant connection on the
 tangent bundle can be thought of as an $H$--equivariant extension $A:\;\g
 \,\longrightarrow\,\End\,\g/\h$ of the adjoint representation $\star$ of
 $\h$ on $\g/\h$. We need to show that $[\,N,\,X\,]\,\in\,\n$ for all $N\,
 \in\,\n$ and $X\,\in\,\g$ or equivalently $[\,[\,N,\,X\,],\,Y\,]\,\equiv
 \,0\;\mod\;\h$ for all $X,\,Y\,\in\,\g$, because $[\,N,\,X\,]\,\in\,\h$
 by the very definition of $\n$. Since the left invariant connection $A$
 extends the infinitesimal representation
 $$
  [\,[\,N,\,X\,],\,Y\,]
  \;\;\equiv\;\;
  A_{[\,N,\,X\,]}(\,Y\,+\,\h\,)
  \;\;\equiv\;\;
  [\;N,\;A_X(\,Y\,+\,\h\,)\;]\,-\,A_X(\,[\,N,\,Y\,]\,+\,\h\,)
 $$
 modulo $\h$, where the second congruence is the infinitesimal version of
 $H$--equivariance for the left invariant connection $A$ under $N\,\in\,\h$.
 Evidently the right hand side vanishes, because $[\,N,\,Y\,]\,\in\,\h$ and
 $[\,N,\,A_X(Y+\h)\,]\,\in\,\h$ by assumption.
 \qed

 \pfill
 In general the kernel $\n$ of the adjoint representation of $\h$ on
 $\g/\h$ fails to be an ideal of $\g$ and thus obstructs the existence
 of a left invariant connection on the tangent bundle and in turn a
 left invariant principal connection on a homogeneous space $G/H$.
 Homogeneous spaces $G/H$ admitting a left invariant connection on their
 tangent bundle are really rather special in the class of all homogeneous
 spaces, because $G$ will not be a subgroup of some affine group in general.
 Of course this simple observation does not preclude the existence of left
 invariant connections on other homogeneous vector bundles. The infinitesimal
 representation $\star:\,\g\longrightarrow\End\,\Sigma$ of a representation
 $\Sigma$ of $G$ for example is always a left invariant connection on a
 homogeneous vector bundle over $G/H$ with base point fiber $\Sigma_{eH}(G/H)
 \,=\,\Sigma$.

 \begin{Corollary}[Kernel of the Adjoint Representation]
 \hfill\label{eff}\break
  If a homogeneous space $G/H$ carries a left invariant connection on its
  tangent bundle, then we may assume without loss of generality that the
  adjoint representation $\star:\,\h\longrightarrow\End\g/\h$ is faithful.
  Namely under this assumption the kernel $\n$ of the adjoint representation
  is the Lie algebra of the maximal normal subgroup $N$ of $G$ contained in
  $H$ so that we may present $G/H$ alternatively as $(G/N)\,/\,(H/N)$ with
  faithful adjoint representation $\star:\,\h/\n\longrightarrow\End\,\g/\h$.
 \end{Corollary}

 \pfill
 Recall now that the group $G$ acts on the vector space $\Gamma(\,\Sigma
 (G/H)\,)$ of a homogeneous vector bundle via $(L_\gamma s)(gH)\,:=\,\gamma
 \,\star\,s(\gamma^{-1}gH)$, sections invariant under this representation
 in the sense $L_\gamma s\,=\,s$ for all $\gamma\,\in\,G$ are called left
 invariant sections. The value of such a left invariant section $s$ in the
 base point $eH\,\in\,G/H$ is necessarily an invariant vector in the base
 point fiber representation $\Sigma\,=\,\Sigma_{eH}(G/H)$, since for
 $h\,\in\,H$ it is necessarily true that:
 $$
  s(\,eH\,)\;\;=\;\;(\,L_hs\,)(\,eH\,)\;\;=\;\;h\,\star\,s(\,h^{-1}H\,)
 $$
 This argument is the central observation leading to a classification
 of left invariant sections:

 \begin{Lemma}[Characterization of Left Invariant Sections]
 \hfill\label{cli}\break
  Evaluation at the base point provides an isomorphism between the vector
  space $\Gamma(\,\Sigma(G/H)\,)^G$ of left invariant sections of a
  homogeneous vector bundle $\Sigma(G/H)$ on a homogeneous space $G/H$
  and the subspace $\Sigma^H$ of $H$--invariant vectors in the
  $H$--representation $\Sigma\,:=\,\Sigma_{eH}(G/H)$:
  $$
   \mathrm{ev}_{eH}:\quad[\;\Gamma(\;\Sigma(G/H)\;)\;]^G
   \;\longrightarrow\;[\;\Sigma\;]^H,\qquad s\;\longmapsto\;s(\,eH\,)
  $$
  The inverse isomorphism associates to $s_{eH}\,\in\,\Sigma^H$ the well
  defined section $s(\,gH\,)\,:=\,g\,\star\,s_{eH}$.
 \end{Lemma}

 \pfill
 With the characterization of left invariant sections at hand we can
 eventually make the relation between the algebraic and the analytic
 definition of left invariant connections somewhat more precise. Observing
 that the jet bundles $\Jet^r\Sigma(G/H)$ of all orders $r\,\geq\,0$
 of a homogeneous vector bundle $\Sigma(G/H)$ are naturally homogeneous
 under the left $G$--action
 $$
  \gamma\,\star\,\mathrm{jet}^r_{gH}s
  \;\;:=\;\;
  \mathrm{jet}^r_{\gamma gH}(\,L_\gamma s\,)
  \;\;=\;\;
  \mathrm{jet}^r_{\gamma gH}\Big[\;\gamma\tilde gH\,\longmapsto\,
  \gamma\,\star\,s(\,\tilde g\,H\,)\;\Big]
 $$
 we can convert every left invariant differential operator $D:\,\Gamma
 (\,\Sigma(G/H)\,)\longrightarrow\Gamma(\,\tilde\Sigma(G/H)\,)$ of order
 $r\,\geq\,0$ between sections of homogeneous vector bundles $\Sigma(G/H)$
 and $\tilde\Sigma(G/H)$ into a left invariant section of the homogeneous
 vector bundle $\Hom(\,\Jet^r\Sigma(G/H),\,\tilde\Sigma(G/H)\,)$. According
 to Lemma \ref{cli} this left invariant total symbol section corresponds to
 an $H$--equivariant homomorphism $\Jet^r\Sigma\longrightarrow\tilde\Sigma$
 between the respective representations $\Jet^r\Sigma$ and $\tilde\Sigma$
 of $H$. In order to classify left invariant connections using this idea
 we need a suitable model, say
 $$
  \Jet^1\Sigma
  \;\;:=\;\;
  \ker\Big(\;(\,\g^*\otimes\Sigma\,)\,\oplus\,\Sigma
  \;\longrightarrow\;\h^*\,\otimes\,\Sigma,\quad ds\,\oplus\,s
  \;\longmapsto\;\mathrm{res}_\h(\,ds\,)\,+\,\delta s\;\Big)
 $$
 for the representation $\Jet^1\Sigma$, where $(\delta s)(Z)\,:=\,Z\star s$
 for all $Z\,\in\,\h$. This model fits nicely into the short exact sequence
 of representations to be expected from the symbol sequence
 $$
  0\;\longrightarrow\;(\g/\h)^*\,\otimes\,\Sigma
  \;\stackrel\iota\longrightarrow\;\Jet^1\Sigma
  \;\stackrel\pr\longrightarrow\;\Sigma\;\longrightarrow\;0
 $$
 with $\iota(\,\sigma\,)\,:=\,\sigma\oplus0$ and $\pr(\,ds\oplus s\,)\,:=\,s$.
 In consequence the algebraic Definition \ref{lconn} of left invariant
 connections on homogeneous vector bundles is based on the following
 bijection
 $$
  \nabla(\,ds\,\oplus\,s\,)\;\;=\;\;ds\;+\;A\,s
  \qquad\Longleftrightarrow\qquad
  A\,s\;\;=\;\;\nabla(\,ds\,\oplus\,s\,)\;-\;ds
 $$
 between the $H$--equivariant sections $\nabla:\,\Jet^1\Sigma\longrightarrow
 (\g/\h)^*\otimes\Sigma$ of the symbol sequence and the $H$--equivariant
 extensions $A:\,\g\longrightarrow\End\,\Sigma$ of the infinitesimal
 representation of $\h$ on $\Sigma$.

 \begin{Lemma}[Curvature and Torsion of a Left Invariant Connection]
 \hfill\label{ctlic}\break
  The curvature of a left invariant connection $A:\,\g\longrightarrow
  \End\,\Sigma$ on a homogeneous vector bundle $\Sigma(G/H)$ is the
  $\End\,\Sigma$--valued 2--form $R\,\in\,[\,\L^2(\g/\h)^*\otimes\End
  \,\Sigma\,]^H$ on $\g/\h$ defined by
  $$
   R_{X\,+\,\h,\,Y\,+\,\h}\;\;:=\;\;[\;A_X,\;A_Y\;]\;-\;A_{[\,X,\,Y\,]}
  $$
  for representatives $X,\,Y\,\in\,\g$. In the same vein the torsion of
  a left invariant connection $A:\,\g\longrightarrow\End\,\g/\h$ on the
  tangent bundle $T(G/H)$ of the homogeneous space $G/H$ is defined as the
  $\g/\h$--valued $2$--form $T\,\in\,[\,\L^2(\g/\h)^*\otimes(\g/\h)\,]^H$
  given on representatives $X,\,Y\,\in\,\g$ by:
  $$
   T(\;X\,+\,\h,\,Y\,+\,\h\;)
   \;\;\equiv\;\;
   A_X(\,Y\,+\,\h\,)\;-\;A_Y(\,X\,+\,\h\,)\;-\;[\,X,\,Y\,]
   \qquad\mod\;\;\h
  $$
 \end{Lemma}

 \pfill
 The usual projection from the space of linear connections on the tangent
 bundle $TM$ of a manifold $M$ to the space of torsion free connections
 works just as well for left invariant connections on the tangent bundle
 of a homogeneous space $G/H$. Because the algebraic torsion $T\,\in\,
 [\,\L^2(\g/\h)^*\otimes(\g/\h)\,]^H$ defined in Lemma \ref{ctlic} is
 $H$--invariant, we may use it to modify the left invariant connection
 $A:\,\g\longrightarrow\End\,\g/\h$ to the $H$--equivariant linear map
 \begin{equation}\label{tfpr}
  A^\tf_X(\;Y\,+\,\h\;)
  \;\;:=\;\;
  A_X(\;Y\,+\,\h\;)\;-\;\tfrac12\;T(\;X\,+\,\h,\;Y\,+\,\h\;)
 \end{equation}
 which still extends the infinitesimal isotropy representation $\star$ of
 $\h$ on $\g/\h$. Clearly the torsion $T^\tf$ of this new left invariant
 connection $A^\tf:\,\g\longrightarrow\End\,\g/\h$ vanishes by construction.
 Thinking of the torsion $T(\,X+\h,\,Y+\h\,)\,=:\,T_{X+\h}(\,Y+\h\,)$ as an
 endomorphism valued $1$--form on $\g/\h$ and using $(A_X\star T)_{Y+\h}
 \,:=\,[\,A_X,\,T_{Y+\h}\,]\,-\,T_{A_X(Y+\h)}$ we find for the curvature:
 \begin{eqnarray*}
  \lefteqn{R^\tf_{X\,+\,\h,\;Y\,+\,\h}}\;\;
  &&
  \\[4pt]
  &=&
  [\,A_X-\tfrac12\,T_{X\,+\,\h},\,A_Y-\tfrac12\,T_{Y\,+\,\h}\,]\,-\,
  A_{[X,Y]}\,+\,\tfrac12\,T_{A_X(Y+\h)-A_Y(X+\h)-T(X+\h,Y+\h)}
  \\
  &=&
  R_{X+\h,\,Y+\h}
  \,-\,\tfrac12\Big(\,(A_X\star T)_{Y+\h}\,-\,(A_Y\star T)_{X+\h}\,\Big)
  \,+\,\tfrac14\Big(\,[T_{X+\h},T_{Y+\h}]\,-\,2\,T_{T(X+\h,Y+\h)}\,\Big)
 \end{eqnarray*}
 This formula is remarkably similar to the standard formula for the
 curvature $R^\tf$ of the torsion free projection $\nabla^\tf\,:=\,
 \nabla-\frac12T$ of a connection on the tangent bundle a manifold $M$
 $$
  R^\tf_{X,Y}
  \;\;=\;\;
  R_{X,Y}
  \,-\,\tfrac12\Big(\,(\nabla_XT)_Y\;-\;(\nabla_YT)_X\,\Big)
  \,+\,\tfrac14\Big(\,[T_X,T_Y]\,-\,2\,T_{T(X,Y)}\,\Big)
 $$
 provided we identify $(A_X\star T)_{Y+\h}$ with $(\nabla_XT)_Y$. Generalizing
 this identification we find:

 \begin{Remark}[Covariant Derivatives of Left Invariant Sections]
 \hfill\label{icd}\break
  The covariant derivative of a left invariant section $s\,\in\,\Gamma(\,
  \Sigma(G/H)\,)^G$ of a homogeneous vector bundle $\Sigma(G/H)$ under a
  left invariant connection $A:\,\g\longrightarrow\End\,\Sigma$ is a left
  invariant section $\nabla^As$ of the homogeneous vector bundle $T^*(G/H)
  \otimes\Sigma(G/H)$. In particular the identification $\Gamma(\,\Sigma
  (G/H)\,)^G\,\cong\,[\Sigma]^H$ turns the left invariant connection into
  the linear map
  $$
   A\,\bm\,:
   \quad[\;\Sigma\;]^H\;\longrightarrow\;[\;(\g/\h)^*\,\otimes\,\Sigma\;]^H,
   \qquad s\;\longmapsto\;A\,\bm\,s
  $$
  well--defined by $(A\,\bm\,s)_{X+\h}\,:=\,A_Xs$ for $X\,\in\,\g$, because
  $A_Xs\,=\,X\star s\,=\,0$ for all $X\,\in\,\h$. 
 \end{Remark}

 \pfill
 A pleasant aspect of this description of the covariant derivative of left
 invariant sections is that it can be iterated with ease provided we we fix
 an auxiliary left invariant section $A^\aux:\,\g\longrightarrow\End\,\g/\h$
 on the tangent bundle $T(G/H)$ of the homogeneous space $G/H$. The standard
 skew symmetrization for iterated covariant derivatives for example becomes
 $$
  \mathrm{alt}_{12}\Big(\;(A,A^\aux)\,\bm\,(\;A\,\bm\,s\;)\;\Big)
  \;\;=\;\;
  R\,\bm\,s\;-\;A_{T^\aux}\,\bm\,s
 $$
 where $\mathrm{alt}_{12}:\,V^*\otimes(V^*\otimes\Sigma)\longrightarrow
 (V^*\otimes V^*)\otimes\Sigma$ is the skew symmetrization in the first
 two arguments and $T^\aux$ the torsion of the auxiliary connection $A^\aux$.
\section{Formal Affine Homogeneous Spaces}
\label{afs}
 Abstracting the concept of affine homogeneous spaces developped so far
 into a purely algebraic concept is only possibly, if we agree to ignore
 the global aspects of a homogeneous space $G/H$ considered as a manifold
 and focus on the pair $\g\,\supset\,\h$ of Lie algebras instead. Since
 it is no longer feasible in this algebraic context to define a left
 invariant connection as an $H$--equivariant linear map, we are lead
 to define a formal affine homogeneous space as a pair $\g\,\supset\,\h$
 of Lie algebras endowed with an $\h$--equivariant extension $A:\,\g
 \longrightarrow\End\,\g/\h$ of the infinitesimal isotropy representation
 $\star:\,\h\longrightarrow\End\,\g/\h$. In this section we will associate
 two isomorphic skew algebras $(\End\,\g/\h)\,\oplus_{\h,A}\g\,\cong\,
 (\End\,\g/\h)\,\oplus_{R,T}(\g/\h)$ to such a formal affine homogeneous
 space with the characteristic property that the quotient $\g/\n$ of $\g$
 by the maximal ideal $\n\,\subset\,\h$ contained in $\h$ becomes a Lie
 subalgebra of these skew algebras.

 \begin{Definition}[Formal Affine Homogeneous Spaces]
 \hfill\label{fahs}\break
  A formal affine homogeneous space is a pair $\g\,\supset\,\h$ of
  Lie algebras endowed with a formal left invariant connection on
  its isotropy representation, this is an $\h$--equivariant linear
  map $A:\,\g\longrightarrow\End\,\g/\h$ extending the isotropy
  representation $\star:\,\h\longrightarrow\End\,\g/\h$. The curvature
  $R\,\in\,[\L^2(\g/\h)^*\otimes(\End\,\g/\h)]^\h$ and the torsion
  $T\,\in\,[\L^2(\g/\h)^*\otimes(\g/\h)]^\h$ of a formal affine
  homogeneous space are defined by the same formulas as for an
  actual homogeneous space:
  \begin{eqnarray*}
   R_{X\,+\,\h,\,Y\,+\,\h}\qquad\quad
   &:=&
   [\;A_X,\;A_Y\;]\;-\;A_{[\,X,\,Y\,]}
   \\[2pt]
   T(\,X+\h,\,Y+\h\,)
   &:=&
   A_X(\,Y\,+\,\h\,)\;-\;A_Y(\,X\,+\,\h\,)\;-\;[\,X,\,Y\,]\;+\;\h
  \end{eqnarray*}
 \end{Definition}

 \pfill
 Evidently every affine homogeneous space $G/H$ endowed with a left
 invariant connection $A:\,\g\longrightarrow\End\,\g/\h$ on its tangent
 bundle $T(G/H)$ defines a formal affine homogeneous spaces $\g\,\supset
 \,\h$ with formal connection $A$. On the other hand we may encounter
 serious problems in integrating a formal affine homogeneous space to
 an actual homogeneous space, because the subset $\exp\,\h\,\subset\,G$
 of the simply connected Lie group $G$ with Lie algebra $\g$ does not
 in general generate a closed subgroup $H\,\subset\,G$, quite simple
 and beautiful conterexamples in this direction have been constructed
 by Kowalski \cite{kov}.

 In the case $H$ fails to be a closed subgroup of $G$ the formal affine
 homogeneous space $\g\,\supset\,\h$ does not integrate to a true affine
 homogeneous space $G/H$. In order to define some kind of surrogate we
 may consider the closure $\overline{H}\,\supset\,H$ of $H$ in $G$, which
 is a Lie subgroup of $G$ with its own Lie algebra $\overline{\h}\,\supset
 \,\h$. Since the adjoint action $\star:\,G\times\g\longrightarrow\g$ of
 the Lie group $G$ on its Lie algebra $\g$ is continuous, the $H$--invariant
 subspace $\h\,\subset\,\g$ is actually $\overline{H}$--invariant so that
 $\overline\h/\h$ is a Lie algebra. The resulting short exact sequence of
 representations of $\overline{H}$
 $$
  0
  \;\longrightarrow\;\overline\h/_{\displaystyle\h}
  \;\longrightarrow\;\g/_{\displaystyle\h}
  \;\longrightarrow\;\g/_{\displaystyle\overline\h}
  \;\longrightarrow\;0
 $$
 corresponds to a short exact sequence of homogeneous vector bundles on
 $G/\overline{H}$ involving the tangent bundle $T(G/\overline{H})$ modelled
 on $\g/\overline\h$ and the homogeneous Lie algebra bundle modelled on
 $\overline\h/\h$. Hence the homogeneous vector bundle on $G/\overline{H}$
 modelled on $\g/\h$ is a transitive Lie algebroid bundle ``$T(G/H)$''
 endowed with a left invariant connection $\nabla$, whose curvature $R$
 and torsion $T$ equal the formal curvature and torsion defined above.
 In this context we recall that the torsion is actually defined for linear
 connections on transitive Lie algebroids.

 \begin{Definition}[Skew Algebra associated to Connection]
 \label{ghalg}\hfill\break
  Consider a left invariant, not necessarily torsion free connection
  $A:\,\g\longrightarrow\End\,\g/\h$ on the isotropy representation
  $\g/\h$ of a pair $\g\,\supset\,\h$ of Lie algebras. The quotient
  of $(\End\,\g/\h)\oplus\g$
  $$
   (\End\,\g/\h)\;\oplus_{\h,A}\,\g
   \;\;:=\;\;
   (\End\,\g/\h)\;\oplus\,\g\;/_{\displaystyle\{\;\;(\,-\,H\,\star\,)
   \,\oplus\,H\;\;|\;\;H\,\in\,\h\;\;\}}
  $$
  by the diagonal subspace $\h$ can be endowed with a skew symmetric
  bilinear bracket via:
  \begin{eqnarray*}
   [\,\X\oplus_{\h,A}X,\,\Y\oplus_{\h,A}Y\,]
   &:=&
   [\,\X,\,\Y\,]\oplus_{\h,A}[\,X,\,Y\,]
   \;+\;(\,[\;\X,\;A_Y\;]\,-\,A_{\X Y}\,)\oplus_{\h,A}\X Y
   \\
   &&
   \qquad\qquad\qquad\qquad\;\;
   \;\,-\,\;(\,[\,\Y,\,A_X\,]\,-\,A_{\Y X}\,)\oplus_{\h,A}\Y X
  \end{eqnarray*}
  The notation $(\End\,\g/\h)\oplus_{\h,\,A}\g$ reflects the dependence
  of the resulting skew algebra on $A$.
 \end{Definition}

 \noindent
 Some thoughts should be spent on the interpretation of the terms $\X Y$
 and $\Y X$ in the definition of the bracket above, which are used as if
 the classes $\X(\,Y\,+\,\h\,)$ and $\Y(\,X\,+\,\h\,)$ in $\g/\h$ were well
 defined elements of $\g$. Nevertheless the resulting ambiguities cancel
 out in the quotient $(\End\,\g/\h)\oplus_{\h,A}\g$ of $(\End\,\g/\h)
 \oplus\g$ by the diagonal $\{\;(-H\star\,)\oplus H\;|\;\;H\,\in\,\h\;\}$
 $$
  (\,[\,\X,\,A_Y\,]\,-\,A_{\X Y}\,)\oplus_{\h,A}\X Y
  \;\;=\;\;[\,\X,\,A_Y\,]\oplus_{\h,A}0\;+\;(\,-\,A_{\X Y}\,)\oplus_{\h,A}\X Y
 $$
 as long as we take the same representatives $\X Y$ and $\Y X$ in $\g$ for
 the classes $\X(Y+\h)$ and $\Y(X+\h)$ in $\g/\h$ on both sides of $\oplus$.
 With this proviso the bracket of two elements in $(\End\,\g/\h)\oplus\g$
 is well defined in the quotient $(\End\,\g/\h)\oplus_{\h,A}\g$. The bracket
 will descend to a skew algebra structure on $(\End\,\g/\h)\oplus_{\h,A}\g$,
 if all elements representing $0$ in the quotient have vanishing brackets
 with all other elements of $(\End\,\g/\h)\oplus\g$. Observing that the
 classes $H\star(Y\,+\,\h)$ and $\Y(H\,+\,\h)$ for given $H\,\in\,\h$ are
 represented in $\g$ by $[H,Y]$ and $0$ we find
 \begin{eqnarray*}
  \lefteqn{[\;(-H\star)\,\oplus\,H,\;\Y\,\oplus\,Y\;]}
  &&
  \\[2pt]
  &=&
  [\,-\,H\star,\,\Y\,]\oplus_{\h,A}[\,H,\,Y\,]
  \;+\;\Big(\,[\,-\,H\star,\,A_Y\,]\;+\;A_{[H,Y]}
  \;-\;[\,\Y,\,A_H\,]\,\Big)\oplus_{\h,A}\Big(\,-\,[\,H,\,Y\,]\,\Big)
  \\
  &=&-\;\Big(\;[\;H\,\star,\;A_Y\;]\;-\;A_{[H,Y]}\;\Big)\oplus_{\h,A}0
 \end{eqnarray*}
 using $A_H\,=\,H\star$ in the second line. Due to the characteristic
 infinitesimal $\h$--equivariance of the left invariant connection
 $A:\,\g\longrightarrow\End\,\g/\h$ the right hand side vanishes so
 that the bracket is in fact well--defined on the quotient
 $(\End\,\g/\h)\oplus_{\h,A}\g$.

 The inclusion of the first summand $\X\,\longmapsto\,\X\oplus_{\h,A}0$ is
 evidently an injective skew algebra homomorphism turning $\End\,\g/\h$ into
 a Lie subalgebra of $(\End\,\g/\h)\oplus_{\h,A}\g$. Things are slightly more
 complicated for the inclusion $X\longmapsto 0\oplus_{\h,A}X$ of the second
 summand, which may well fail to be injective. Its kernel however agrees with
 the kernel of the isotropy representation $\star:\,\h\longrightarrow\End\,
 \g/\h,\,H\longmapsto H\,\star,$ which in turn agrees with the maximal ideal
 $\n\,\subset\,\h$ of $\g$ contained in $\h$ according to Lemma \ref{eff}.
 In this way the skew algebra $(\End\,\g/\h)\,\oplus_{\h,A}\g$ is generated
 as a vector space by its two Lie subalgebras $\End\,\g/\h$ and $\g/\n$,
 nevertheless the Jacobi identity will not hold true for arbitrary elements
 of $(\End\,\g/\h)\,\oplus_{\h,A}\g$.

 \begin{Definition}[Skew Algebra associated to Curvature--Torsion]
 \label{rtalg}\hfill\break
  Consider a left invariant connection $A:\,\g\longrightarrow\End\,\g/\h$
  on the isotropy representation $\g/\h$ of a pair $\g\,\supset\,\h$ of
  Lie algebras with associated curvature $R\,\in\,[\,\L^2(\g/\h)^*\otimes
  \End\,\g/\h\,]^\h$ and torsion $T\,\in\,[\,\L^2(\g/\h)^*\otimes\g/\h\,]^\h$.
  The direct sum $(\End\,\g/\h)\oplus(\g/\h)$ of vector spaces is actually a
  skew algebra denoted by $(\End\,\g/\h)\oplus_{R,T}(\g/\h)$ under the
  following skew bilinear bracket: 
  $$
   [\;\X\,\oplus_{R,T}x,\;\Y\,\oplus_{R,T}y\;]
   \;\;:=\;\;
   \Big(\;[\,\X,\,\Y\,]\;-\;R_{x,\,y}\;\Big)
   \,\oplus_{R,T}\Big(\;\X y\;-\;\Y x\;-\;T(\,x,\,y\,)\;\Big)
  $$
 \end{Definition}

 \noindent
 Up to the additional terms $-\,R_{x,\,y}$ and $-\,T(x,y)$ the skew bracket
 on $(\End\,\g/\h)\oplus_{R,T}(\g/\h)$ agrees with the Lie bracket on the
 semidirect product $(\End\,\g/\h)\,\oplus\,(\g/\h)$ of the Lie algebra
 $\End\,\g/\h$ with its representation $\g/\h$. The possible failure of
 the Jacobi identity for the bracket on $(\End\,\g/\h)\oplus_{R,T}(\g/\h)$
 is thus due to these two additional terms, for the moment however we postpone
 a detailed analysis of this problem to Section \ref{class}.

 \begin{Lemma}[Skew Algebra Isomorphism]
 \label{isoalg}\hfill\break
  Let $A:\,\g\longrightarrow\End\,(\g/\h)$ be a left invariant connection
  on the isotropy representation $\g/\h$ of a pair $\g\,\supset\,\h$ of Lie
  algebras with curvature $R\,\in\,[\,\L^2(\g/\h)^*\otimes(\End\,\g/\h)\,]^\h$
  and torsion $T\,\in\,[\,\L^2(\g/\h)^*\otimes(\g/\h)\,]^\h$. The connection
  $A$ defines a canonical isomorphism of skew algebras
  \begin{eqnarray*}
   \Phi_A\,:\quad(\End\,\g/\h)\,\oplus_{\h,A}\g
   &\stackrel\cong\longrightarrow&
   (\End\,\g/\h)\,\oplus_{R,T}(\g/\h)
   \\
   \X\oplus_{\h,\,A}X\quad
   &\longmapsto&
   (\X+A_X)\oplus_{R,T}(X+\h)
  \end{eqnarray*}
  between the skew algebras associated to $A$ and $R,\,T$. Composing $\Phi_A$
  with the skew algebra homomorphism $\g\longrightarrow(\End\,\g/\h)\,
  \oplus_{\h,\,A}\g,\,X\longmapsto 0\,\oplus_{\h,\,A}X,$ we may identify
  the quotient $\g/\n$ of $\g$ by the maximal ideal $\n\,\subset\,\g$
  contained in $\h$ with a Lie subalgebra of $(\End\,\g/\h)\,\oplus_{R,\,T}
  (\g/\h)$:
  $$
   \g/\n\;\stackrel\subset\longrightarrow\;(\End\,\g/\h)\,\oplus_{R,T}(\g/\h),
   \qquad X\,+\,\n\;\longmapsto\;A_X\,\oplus_{R,T}(X+\h)
  $$
 \end{Lemma}

 \proof
 A short inspection shows that $\Phi_A$ is well defined on the quotient
 $(\End\,\g/\h)\oplus_{\h,A}\g$ due to $A_H\,=\,H\star$ for all $H\,
 \in\,\h$. Moreover $\Phi_A$ is a linear isomorphism with well--defined
 inverse $\Phi_A^{-1}:\,(\End\,\g/\h)\oplus_{R,T}(\g/\h)\longrightarrow
 (\End\,\g/\h)\oplus_{\h,A}\g$ given explicitly by
 $$
  \Phi_A^{-1}\Big(\;\X\,\oplus_{R,T}x\;\Big)
  \;\;:=\;\;
  (\;\X\;-\;A_X\;)\,\oplus_{\h,A}X
 $$
 where $X\,\in\,\g$ represents the class $x\,\in\,\g/\h$. Eventually $\Phi_A$
 is a homomorphism of algebras
 \begin{eqnarray*}
  \lefteqn{\Big[\;\Phi_A(\;\X\oplus_{\h,A}X\;),\;
   \Phi_A(\;\Y\oplus_{\h,A}Y\;)\;\Big]}\quad
   &&
   \\[2pt]
   &=&
   \Big(\;[\,\X,\,\Y\,]\;+\;[\,A_X,\,\Y\,]\;+\;
    [\,\X,\,A_Y\,]\;+\;[\,A_X,\,A_Y\,]\;-\;R_{X+\h,Y+\h}\;\Big)
   \\[-1pt]
   &&
   \qquad\oplus_{R,T}\Big(\;(\,\X+A_X\,)(\,Y+\h\,)
   \;-\;(\,\Y+A_Y\,)(\,X+\h\,)\;-\;T(X+\h,Y+\h)\;\Big)
   \\
   &=&
   \Big(\,[\,\X,\,\Y\,]\,+\,[\,\X,\,A_Y\,]\,-\,[\,\Y,\,A_X\,]
   \,+\,A_{[X,Y]}\,\Big)\oplus_{R,T}\Big(\,\X\,Y\,-\,\Y\,X\,-\,[X,Y]
   \,+\,\h\,\Big)
   \\
   &=&
   \Phi_A\Big(\,(\,[\X,\Y]+[\X,A_Y]-A_{\X\,Y}-[\Y,A_X]+A_{\Y\,X}\,)
   \oplus_{\h,A}(\,\X\,Y-\Y\,X+[X,Y]\,)\;\Big)
 \end{eqnarray*}
 where $\X\,Y,\,\Y\,X\,\in\,\g$ denote fixed representatives for the
 corresponding classes in $\g/\h$.
 \qed
\section{A Parametrization of Formal Affine Spaces}
\label{class}
 In this section we change our point of view completely away from a fixed
 formal affine homogeneous space towards a parametrization of such spaces
 by connection--curvature--torsion triples. In order to compare different
 affine homogeneous spaces we choose a linear isomorphism or frame $F:\,
 V\longrightarrow\g/\h$ to pull back curvature and torsion to $V$ and choose
 a section $\rep:\,\g/\h\longrightarrow\g$ to complement the information
 contained in $R\,\in\,\L^2V^*\otimes\End\,V$ and $T\,\in\,\L^2V^*\otimes V$
 by an additional $A\,\in\,V^*\otimes\End\,V$ describing the connection. A
 generic triple $(A,R,T)$ of this form will certainly not come from a formal
 affine homogeneous space $\g\,\supset\,\h$, the algebraic equations
 characterizing connection--curvature--torsion triples $(A,R,T)$ of formal
 affine homogeneous spaces will be made explicit at the end of this section.

 \pfill
 For the time being let us fix a finite dimensional vector space $V$. 
 Augmenting a formal affine homogeneous space $\g\,\supset\,\h$ of
 dimension $\dim\,\g/\h\,=\,\dim\,V$ with a linear isomorphism or frame
 $F:\,V\longrightarrow\g/\h$ allows us to think of curvature and torsion
 of the left invariant connection $A:\,\g\longrightarrow\End\,\g/\h$ as
 elements $R\,\in\,\L^2V^*\otimes\End\,V$ and $T\,\in\,\L^2V^*\otimes V$.
 Nevertheless the resulting curvature--torsion tuple $(\,R,\,T\,)$ does
 not describe the original formal affine homogeneous space completely.
 For this reason we choose in addition a section $\rep: \,\g/\h
 \longrightarrow\g$ of the canonical projection to $\g/\h$ in order to
 capture the information contained in the connection $A$ in a linear
 map $A\circ\rep:\,\g/\h\longrightarrow\End\,\g/\h$, which becomes
 under $F$ the $\End\,V$--valued $1$--form $A\,\in\,V^*\otimes\End\,V$
 on $V$ still called connection:

 \begin{Definition}[Connection--Curvature--Torsion Triples]
 \hfill\label{cct}\break
  A connection--curvature--torsion triple on a vector space $V$ is a triple
  of the form:
  $$
   (\;A,\;R,\;T\;)\;\;\in\;\;
   (\;V^*\,\otimes\,\End\;V\;)
   \;\times\;
   (\;\L^2V^*\,\otimes\,\End\;V\;)
   \;\times\;
   (\;\L^2V^*\,\otimes\,V\;)
  $$
  Every such a triple $(A,R,T)$ endows $\End\,V\oplus_{R,T}V\,:=\,\End\,V
  \oplus V$ with the skew bracket:
  $$
   [\;\X\,\oplus_{R,T}x,\;\Y\,\oplus_{R,T}y\;]
   \;\;:=\;\;
   \Big(\;[\;\X,\;\Y\;]\;-\;R_{x,y}\;\Big)\,\oplus_{R,T}
   \Big(\;\X\,y\;-\;\Y\,x\;-\;T(\,x,\,y\,)\;\Big)
  $$
 \end{Definition}

 \pfill
 The connection--curvature--torsion triples $(A,R,T)$ coming from actual
 formal affine homogeneous spaces $\g\,\supset\,\h$ augmented by frames
 $F:\,V\longrightarrow\g/\h$ and sections $\rep:\,\g/\h\longrightarrow\g$
 are characterized by the fact that the skew algebra $\End\,V\oplus_{R,T}V$
 contains a Lie subalgebra
 $$
  \g/\n\;\stackrel\subset\longrightarrow\;(\End\,\g/\h)\,\oplus_{R,T}(\g/\h)
  \;\stackrel\cong\longrightarrow\;\End\,V\,\oplus_{R,T}V
 $$
 isomorphic to $\g/\n$ according to Lemma \ref{isoalg}, which contains
 the image of the extension of $A$
 $$
  A^{\mathrm{ext}}:\quad V\;\longrightarrow\;\End\,V\,\oplus_{R,T}V,
  \qquad x\;\longmapsto\;A_x\,\oplus_{R,T}x
 $$
 and thus projects onto $V$ under the canonical projection $\End\,V
 \oplus_{R,T}V\longrightarrow V$. Conversely:

 \begin{Definition}[Isotropy Algebra and Tautological Connection]
 \label{taut}\hfill\break
  Consider a Lie subalgebra $\g$ of the skew algebra $\End\,V\oplus_{R,T}V$
  associated to a curvature--torsion tuple $(R,T)\,\in\,\L^2V^*\otimes\End\,V
  \times\,\L^2V^*\otimes V$, which projects surjectively onto $V$ under:
  $$
   \End\,V\,\oplus_{R,T}V\;\longrightarrow\;V,
   \qquad\X\,\oplus_{R,T}x\;\longmapsto\;x
  $$
  Defining the isotropy algebra $\h\,:=\,\g\,\cap\,\End\,V$ as the kernel
  of this projection we observe that the canonical projection becomes an
  isomorphism $\g/\h\stackrel\cong\longrightarrow V$ of the isotropy
  representation of the pair $\g\,\supset\,\h$ with $V$, with this proviso
  the canonical projection to the first summand
  $$
   A^{\mathrm{taut}}:
   \quad\End\,V\,\oplus_{R,T}V\;\longrightarrow\;\End\,V,
   \qquad\X\,\oplus_{R,T}x\;\longmapsto\;\X
  $$
  becomes the tautological left invariant connection $A^{\mathrm{taut}}:\,
  \g\,\longrightarrow\,\End\,V$ on $V\,\cong\,\g/\h$.
 \end{Definition}

 \noindent
 Although or perhaps because all arguments and calculations involving the
 tautological connection $A^{\mathrm{taut}}$ eventually reduce to tautologies
 it is difficult and slightly confusing to use it in explicit calculations.
 Nevertheless the tautological connection $A^{\mathrm{taut}}$ is
 $\h$--equivariant to due $[\X\,\oplus_{R,T}0,\,\Y\,\oplus_{R,T}y]
 \,=\,[\X,\Y]\,\oplus_{R,T}\X\,y$ and its curvature agrees with
 $R\,\in\,\L^2V^*\otimes\End\,V$:
 $$
  [\;A^{\mathrm{taut}}_{\X\,\oplus_{R,T}x},\;
   A^{\mathrm{taut}}_{\Y\,\oplus_{R,T}y}\;]
  \;-\;
   A^{\mathrm{taut}}_{[\,\X\oplus_{R,T}x,\,\Y\oplus_{R,T}y\,]}
  \;\;=\;\;
  [\,\X,\,\Y\,]\;-\;\Big(\;[\,\X,\,\Y\,]\,-\,R_{x,y}\;\Big)
  \;\;=\;\;
  R_{x,y}
 $$
 A very similar argument can be made to verify that the torsion of the
 tautological connection agrees with $T\,\in\,\L^2V^*\otimes V$. In
 consequence every Lie subalgebra $\g\,\subset\,\End\,V\oplus_{R,T}V$
 projecting surjectively onto $V$ defines a formal affine homogeneous
 space $\g\,\supset\,\h$ with tautological frame $F:\,V\longrightarrow
 \g/\h$ under the tautological left invariant connection $A^{\mathrm{taut}}:\,
 \g\longrightarrow\End\,V$. Motivated by Lemma \ref{icd} we now define
 for every representation $\Sigma$ of the Lie algebra $\End\,V$ and every
 $d\,\geq\,0$ the bilinear operation
 $$
  \bm:\quad{\textstyle\bigotimes}^dV^*\,\otimes\,\End\,V
  \,\times\,{\textstyle\bigotimes}^\bullet V^*\,\otimes\,\Sigma
  \;\longrightarrow\;
  {\textstyle\bigotimes}^{\bullet+d}V^*\,\otimes\,\Sigma
 $$
 by setting
 $$
  (\,Q\,\bm\,s\,)(\;x_1,\,\ldots,\,x_d;\,y_1,\,y_2,\,\ldots\;)
  \;\;:=\;\;
  (\,Q_{x_1,\,\ldots,\,x_d}\,\star\,s\,)(\;y_1,\,y_2,\,\ldots\;)
 $$
 where $\star$ denotes the tensor product representation of $\End\,V$
 on $\bigotimes^\bullet V^*\otimes\Sigma$.

 \begin{Definition}[Formal Covariant Derivatives of $R$ and $T$]
 \hfill\label{fid}\break
  For a given connection--curvature--torsion triple $(A,R,T)$ on a
  finite--dimensional vector space $V$ we define the formal iterated
  covariant derivatives of $R$ and $T$ by setting:
  \begin{eqnarray*}
   \nabla^rT
   &:=&
   \underbrace{A\,\bm\,(\,A\,\bm\,(\,\ldots\,(\,A}_{r\;\;\mathrm{times}}
   \,\bm\,T\,)\,\ldots\,)\,)
   \;\;\in\;\;
   {\textstyle\bigotimes}^r\,V^*\,\otimes\,(\,\L^2V^*\,\otimes\,V\,)
   \\
   \nabla^rR
   &:=&
   \underbrace{A\,\bm\,(\,A\,\bm\,(\,\ldots\,(\,A}_{r\;\;\mathrm{times}}
   \,\bm\,R\,)\,\ldots\,)\,)
   \;\;\in\;\;
   {\textstyle\bigotimes}^r\,V^*\,\otimes\,(\,\L^2V^*\,\otimes\,\End\;V\,)
  \end{eqnarray*}
 \end{Definition}

 \begin{Definition}[Stabilizer Filtration]
 \hfill\label{sfilt}\break
  Associated to every connection--curvature torsion triple is the strictly
  descending filtration
  $$
   \End\,V\;=\;\ldots\;=\;\h_{-2}\;=\;\h_{-1}\;\supsetneq\;\h_0
   \;\supsetneq\;\ldots\;\supsetneq\;\h_{s-1}
   \;\supsetneq\;\h_s\;=\;\h_{s+1}\;=\;\ldots\;=\;\h_\infty
  $$
  of $\End\,V$ defined recursively by $\h_0\,:=\,\stab\,R\,\cap\,\stab\,T
  \,\subset\,\End\,V$ and for all $r\,\geq\,0$:
  $$
   \h_{r+1}
   \;\;:=\;\;
   \{\;\;\X\,\in\,\h_r\;\;|\;\;[\,\X,\,A_x\,]\,\equiv\,A_{\X\,x}
   \;\;\mod\;\h_r\textrm{\ \ for all\ }x\,\in\,V\;\;\}
  $$
  The minimal $s\,\geq\,-1$ with equality $\h_s\,=\,\h_{s+1}$ is called the
  Singer invariant of $(A,R,T)$.
 \end{Definition}

 \pfill
 By its recursive definition this filtration is strictly descending in the
 sense that it becomes stationary at the first equality $\h_s\,=\,\h_{s+1}$
 of successive filtration steps for some $s\,\geq\,0$, in passing we observe
 that this argument works even for $s\,=\,-1$, because the equality $\h_{-1}
 \,=\,\h_0$ renders the stipulated congruences $[\X,A_x]\,\equiv\,A_{\X x}$
 modulo $\End\,V$ trivial. A better way to understand the recursive definition
 of the stabilizer filtration is to observe that for every subalgebra $\hat\h
 \,\subset\,\End\,V$ the vector space $V^*\otimes(\End\,V/\hat\h)$ is actually
 a representation of $\hat\h$:
 $$
  (\,\X\,\star\,A\,)_x
  \;\;:=\;\;
  \X\,\star\,A_x\;-\;A_{\X\,\star\,x}\;+\;\hat\h
  \;\;=\;\;
  [\,\X,\,A_x\,]\;-\;A_{\X\,x}\;+\;\hat\h
 $$
 By induction it thus follows that all subspaces $\h_r\,\subset\,\End\,V$ of
 the stabilizer filtration are subalgebras, in fact this is true by definition
 for $\h_0\,:=\,\stab\,R\,\cap\,\stab\,T$, whereas $\h_{r+1}$ for $r\,\geq\,0$
 is the stabilizer subalgebra of the class $A\,+\,V^*\otimes\h_r$ represented
 by $A$ in the $\h_r$--representation $V^*\otimes(\End\,V/\h_r)$. A remarkable
 consequence of this observation is the following very neat interpretation
 of the stabilizer filtration in terms of covariant derivatives:

 \begin{Lemma}[Interpretation of the Stabilizer Filtration]
 \hfill\label{if}\break
  The strictly descending filtration $\h_\bullet$ of $\End\,V$ associated
  to a connection--curvature--torsion triple $(A,R,T)$ can be interpreted
  geometrically as the filtration given by the joint stabilizers
  $$
   \h_r
   \;\;=\;\;
   \stab(\;R\,\oplus\,\nabla R\,\oplus\,\ldots\,\oplus\,\nabla^rR\;)
   \;\cap\;
   \stab(\;T\,\oplus\,\nabla T\,\oplus\,\ldots\,\oplus\,\nabla^rT\;)
  $$
  of the iterated covariant derivatives of curvature $R$ and torsion $T$
  up to order $r\,\geq\,0$.
 \end{Lemma}

 \pfill
 Actually this Lemma is a special case of a more general fact, for every
 representation $\Sigma$ of the Lie algebra $\End\,V$ and every $s\,\in\,
 \Sigma$ with stabilizer $\stab\,s\,\subset\,\End\,V$ the linear map
 $$
  V^*\,\otimes\,(\,\End\,V\,/_{\displaystyle\stab\;s}\,)
  \;\longrightarrow\;V^*\,\otimes\,\Sigma,\qquad
  A\;\longmapsto\;A\,\bm\,s
 $$
 defined by $(A\,\bm\,s)_x\,:=\,A_x\,\star\,s$ is well--defined, injective
 and equivariant under $\stab\,s$. Since the bilinear operation $\bm$ is
 naturally defined and thus equivariant under $\End\,V$ in the sense
 $$
  \X\,\star\,(\,A\,\bm\,s\,)
  \;\;=\;\;
  (\,\X\,\star\,A\,)\,\bm\,s\;+\;A\,\bm\,(\,\X\,\star\,s\,)
 $$
 for $\X\,\in\,\End\,V$ we find that $\X\,\in\,\stab\,s$ stabilizes
 $A\,\in\,V^*\otimes(\End\,V/\stab\,s)$, if and only if:
 $$
  (\;\X\,\star\,A\,)\,\bm\,s\;\;=0\;\;=\;\;\X\,\star\,(\,A\,\bm\,s\,)
 $$
 Using this general argument the statement of the Lemma follows by an easy
 induction based on the definition $\h_0\,:=\,\stab\,R\,\cap\,\stab\,T$ as
 well as on the recursive definition of $\h_{r+1}$ as the stabilizer of the
 connection class $A\,+\,V^*\otimes\h_r\,\in\,V^*\otimes(\End\,V/\h_r)$.

 \pfill
 For a little interlude in our algebraic considerations let us now recall
 the definition of the twisted exterior derivative associated to a connection
 $\nabla$ on a vector bundle $\Sigma M$ over a smooth manifold $M$. The
 twisted exterior derivative $d^\nabla$ is a first order differential
 operator
 $$
  d^\nabla:\quad
  \Gamma(\;\L^\bullet T^*M\,\otimes\,\Sigma M\;)\;\longrightarrow\;
  \Gamma(\;\L^{\bullet+1}T^*M\otimes\Sigma M\;),\qquad
  \omega\;\longmapsto\;d^\nabla\omega
 $$
 on the differential forms on $M$ with values in $\Sigma M$ defined on
 a $\Sigma M$--valued $r$--form $\omega$ by:
 \begin{eqnarray*}
  (\,d^\nabla\omega\,)(\;X_0,\,\ldots,\,X_r\;)
  &:=&
  \sum_{\mu\,=\,0}^r(-1)^\mu\,\nabla_{X_\mu}
  \Big(\;\omega(\;X_0,\,\ldots,\,\widehat{X}_\mu,\,\ldots,\,X_r\;)\;\Big)
  \\
  &&
  \!\!\!-\!\!\!\sum_{0\,\leq\,\mu\,<\,\nu\,\leq\,r}\!\!\!(-1)^{\mu+\nu-1}
  \omega\Big(\;[\,X_\mu,\,X_\nu\,],\,X_0,\,\ldots,\,\widehat{X}_\mu,
  \,\ldots,\,\widehat{X}_\nu,\,\ldots,\,X_r\;\Big)
 \end{eqnarray*}
 A classical formula tells us that with an auxiliary connection $\nabla^\aux$
 on $TM$ we may write
 \begin{eqnarray*}
  d^\nabla\omega(\;X_0,\,\ldots,\,X_r\;)
  &=&
  \sum_{\mu\,=\,0}^r(-1)^\mu\,\Big(\;(\nabla,\nabla^\aux)_{X_\mu}
  \omega\;\Big)(\;X_0,\,\ldots,\,\widehat{X}_\mu,\,\ldots,\,X_r\;)
  \\
  &&
  \!\!\!+\!\!\!
  \sum_{0\,\leq\,\mu\,<\,\nu\,\leq\,r}\!\!\!(-1)^{\mu+\nu-1}\,
  \omega\Big(\,T^\aux(X_\mu,X_\nu),\,X_0,\,\ldots,\,\widehat{X}_\mu,
  \,\ldots,\,\widehat{X}_\nu,\,\ldots,\,X_r\,\Big)
 \end{eqnarray*}
 where $T^\aux\,\in\,\Gamma(\,\L^2T^*M\otimes TM\,)$ is the torsion
 of the auxiliary connection $\nabla^\aux$ and:
 \begin{eqnarray*}
  \lefteqn{\Big(\,(\nabla,\nabla^\aux)_X\omega\,\Big)(\;X_1,\,\ldots,\,X_r\;)}
  \qquad\qquad
  \\[-3pt]
  &:=&
  \nabla_X\Big(\;\omega(\,X_1,\,\ldots,\,X_r\,)\;\Big)
  \;-\;\sum_{\nu\,=\,1}^r\omega\Big(\;
  X_1,\,\ldots,\,\nabla^\aux_XX_\mu,\,\ldots,\,X_r\;\Big)
 \end{eqnarray*}
 In terms of twisted exterior derivatives the First and Second Bianchi
 Identity can be written
 \begin{equation}\label{bi}
  d^{\nabla^\aux}T^\aux\;\;=\;\;R^\aux\wedge\,\id
  \qquad\qquad\qquad
  d^\nabla R\;\;=\;\;0
 \end{equation}
 respectively, where $R$ and $R^\aux$ denote the curvatures of $\nabla$
 and $\nabla^\aux$, whereas:
 $$
  (\,R^\aux\wedge\,\id\,)(\;X,\,Y,\,Z\;)
  \;\;:=\;\;
  R^\aux_{X,\,Y}Z\;+\;R^\aux_{Y,\,Z}X\;+\;R^\aux_{Z,\,X}Y
 $$

 \pfill
 Coming back to formal affine homogeneous spaces and the associated
 connection--curvature--torsion triples $(A,R,T)$ we take the classical
 formula expressing $d^\nabla$ in terms of the connection $(\nabla,
 \nabla^\aux)$ and the torsion $T^\aux$ as a lead to define the twisted
 exterior derivative $d^{(A,T)}$ by
 \begin{eqnarray*}
  (\,d^{(A,T)}\omega\,)(\;x_0,\,\ldots,\,x_r\;)
  &:=&
  \sum_{\mu\,=\,0}^r(-1)^\mu\,
  \Big(\,A_{x_\mu}\star\,\omega\,\Big)(\;x_0,\,\ldots,\,\widehat{x}_\mu,
  \,\ldots,\,x_r\;)
  \\
  &&
  \!\!+\!\!\sum_{0\,\leq\,\mu\,<\,\nu\,\leq\,r}\!\!
  (-1)^{\mu+\nu-1}\omega\Big(\;T(x_\mu,x_\nu),\,x_0,\,\ldots,\,
  \widehat{x}_\mu,\,\ldots,\,\widehat{x}_\nu,\,\ldots,\,x_r\;\Big)
 \end{eqnarray*}
 for every $r$--form $\omega\,\in\,\L^rV^*\otimes\Sigma$ with values
 in a representation $\Sigma$ of the Lie algebra $\End\,V$. Specifically
 for $R\,\in\,\L^2V^*\otimes\End\,V$ the twisted exterior derivative reads
 $$
  (\,d^{(A,T)}R\,)_{x,\,y,\,z}
  \;\;=\;\;
  \Big(\,(A_x\star R)_{y,z}\,+\,R_{T(y,z),\,x}\,\Big)
  \;+\;\textrm{cyclic permutations of\ }x,\,y,\,z
 $$
 and an almost identical formula is valid for $d^{(A,T)}T\,\in\,\L^3V^*
 \otimes V$. Last but not least we define the $V$--valued $3$--form
 $R\,\wedge\,\id\,\in\,\L^3V^*\otimes V$ by $(\,R\wedge\id\,)(x,y,z)
 \,:=\,R_{x,y}z\,+\,R_{y,z}x\,+\,R_{z,x}y$.

 \begin{Definition}[Approximate Curvature]
 \hfill\label{ac}\break
  The approximate curvature tensor of a connection--curvature--torsion
  triple $(A,R,T)$ on a vector space $V$ is defined as an $\End\,V$--valued
  $2$--form $Q(A,T)\,\in\,\L^2V^*\otimes\End\,V$ on $V$ by:
  $$
   Q(\;A,\,T\;)_{x,y}
   \;\;:=\;\;
   [\;A_x,\,A_y\;]\;-\;A_{A_xy\,-\,A_yx\,-\,T(x,y)}
  $$  
 \end{Definition}

 \pfill
 In order to study the failure of the Jacobi identity for the bracket
 of the skew algebra $\End\,V\oplus_{R,T}V$ associated to a
 connection--curvature--torsion triple $(A,R,T)$ we consider
 the trilinear standard Jacobiator defined as an alternating
 $3$--form $\Jac$ on $\End\,V\oplus_{R,T}V$ by:
 \begin{eqnarray*}
  \Jac(\;\X\,\oplus_{R,T}x,\;\Y\,\oplus_{R,T},\;
  \,\mathfrak{Z}\,\oplus_{R,T}z\;)
  &:=&
  +\;[\;\X\,\oplus_{R,T}x,\;
  [\;\Y\,\oplus_{R,T}y,\;\,\mathfrak{Z}\,\oplus_{R,T}z\;]\;]
  \\
  &&
  +\;[\;\Y\,\oplus_{R,T}y,\;
  [\;\,\mathfrak{Z}\,\oplus_{R,T}z,\;\X\,\oplus_{R,T}x\;]\;]
  \\
  &&
  +\;[\;\,\mathfrak{Z}\,\oplus_{R,T}z,\;
  [\;\X\,\oplus_{R,T}x,\;\Y\,\oplus_{R,T}y\;]\;]
 \end{eqnarray*}
 Observing that the bracket with $\X\,\in\,\End\,V$ reproduces the
 infinitesimal representation
 $$
  [\;\X\,\oplus_{R,T}0,\;\Y\,\oplus_{R,T}y\;]
  \;\;=\;\;
  [\;\X,\;\Y\;]\,\oplus_{R,T}\X\,y
  \;\;=\;\;
  \X\,\star\,\Big(\;\Y\,\oplus_{R,T}y\;\Big)
 $$
 of $\End\,V$ on $\End\,V\oplus_{R,T}V$ we may calculate the Jacobiator
 for the special choice
 \begin{eqnarray*}
  \lefteqn{\Jac(\;\X\,\oplus_{R,T}0,\;0\,\oplus_{R,T}y,
   \;0\,\oplus_{R,T}z\;)}\quad
  &&
  \\[3pt]
  &=&
  [\;\X\,\oplus_{R,T}0,\;
  (\,-\,R_{x,\,y}\,)\,\oplus_{R,T}(\,-\,T(\,x,\,y\,)\,)\;]
  \\[2pt]
  &&
  \quad
  \;-\;[\;0\,\oplus_{R,T}y,\;0\,\oplus_{R,T}\X\,z\;]
  \;+\;[\;0\,\oplus_{R,T}z,\;0\,\oplus_{R,T}\X\,y\;]
  \\[2pt]
  &=&
  (\,-\,[\;\X,\,R_{y,\,z}\;]\,+\,R_{\X\,y,\,z}\,+\,R_{y,\,\X\,z}\,)
  \,\oplus_{R,T}(\,-\,\X\,T(\,y,\,z\,)\,+\,T(\,\X\,y,\,z\,)\,+\,
  T(\,y,\,\X\,z\,)\,)
  \\[2pt]
  &=&
  (\,-\,(\,\X\,\star\,R\,)_{y,\,z}\,)\,\oplus_{R,T}
  (\,-\,(\,\X\,\star\,T\,)(\,y,\,z\,)\,)
 \end{eqnarray*}
 of arguments $\X\,\in\,\End\,V$ and $y,\,z\,\in\,V$. Similarly we obtain
 for all three arguments in $V$:
 \begin{eqnarray*}
  \lefteqn{\Jac(\;0\,\oplus_{R,T}x,\;0\,\oplus_{R,T}y,\;
   0\,\oplus_{R,T}z\;)}\quad
  &&
  \\[4pt]
  &=&
  [\;\;0\,\oplus_{R,T}x,
  \;(\,-\,R_{y,\,z}\,)\,\oplus_{R,T}(\,-\,T(\,y,\,z\,)\,)\;\;]
  \;+\;\textrm{cyclic permutations of $x,\,y,\,z$}
  \\[4pt]
  &=&
  R_{x,\,T(\,y,\,z\,)}\,
  \oplus_{R,T}(\;R_{y,\,z}x\,+\,T(\,x,\,T(\,y,\,z\,)\,)\;)
  \;+\;\textrm{cyclic permutations of $x,\,y,\,z$}
  \\[2pt]
  &=&
  \Big(\;-\,d^{(\,0,\,T\,)}R\;\Big)_{x,\,y,\,z}\,\oplus_{R,T}
  \Big(\;R\,\wedge\,\id\;-\;d^{(\,0,\,T\,)}T\;\Big)(\,x,\,y,\,z\,)
 \end{eqnarray*}
 The latter two results feature prominantly in the proof of the following
 lemma:

 \begin{Lemma}[Lie Subalgebras of $\End\,V\oplus_{R,T}V$]
 \hfill\label{sag}\break
  The skew algebra $\End\,V\oplus_{R,T}V$ associated to a
  connection--curvature--torsion triple $(A,R,T)$ on a vector space
  $V$ allows no Lie subalgebra $\g\,\supset\,\{\;A_x\,\oplus_{R,T}x\;|
  \;x\,\in\,V\;\}$ unless the connection--curvature--torsion triple
  $(A,R,T)$ satisfies the First and Second Bianchi Identity:
  $$
   d^{(A,T)}T\;\;=\;\;R\,\wedge\,\id
   \qquad\qquad\qquad
   d^{(A,T)}R\;\;=\;\;0
  $$
  In case the First and Second Bianchi Identity are both satisfied the
  isotropy algebra association $\g\longmapsto\g\,\cap\,\End\,V$ induces
  a bijection between the Lie subalgebras $\g\,\subset\,\End\,V
  \oplus_{R,T}V$ containing $\{\,A_x\oplus_{R,T}x\,|\,x\,\in\,V\,\}$
  and Lie subalgebras $\h\,\subset\,\End\,V$ satisfying
  $$
   \h\;\;\subset\;\;\stab\,R\;\cap\;\stab\,T
   \qquad
   \h\,\star\,A\;\;\subset\;\;V^*\,\otimes\,\h
   \qquad
   Q(\,A,\,T\,)\;\;\equiv\;\;R\;\;\mod\;\L^2V^*\,\otimes\,\h
  $$
  where $\h\,\star\,A\,\subset\,V^*\otimes\h$ is a shorthand for
  $[\,\X,\,A_x\,]\,-\,A_{\X x}\,\equiv\,0$ modulo $\h$ for $x\,\in\,V,
  \,\X\,\in\,\h$.
 \end{Lemma}

 \proof
 Every Lie subalgebra $\g\,\subset\,\End\,V\oplus_{R,T}V$ containing
 $\{\;A_x\,\oplus_{R,T}x\;|\;x\,\in\,V\;\}$ is certainly determined
 by its isotropy algebra $\h\,=\,\g\,\cap\,\End\,V$. Conversely suppose
 that $\h\,\subset\,\End\,V$ is the isotropy algebra of the Lie subalgebra
 $\g$ of $\End\,V\oplus_{R,T}V$ given by:
 \begin{equation}\label{galg}
  \g\;\;:=\;\;
  \h\;+\;\span\,\{\;\;A_x\,\oplus_{R,T}x\;\;|\;\;x\;\in\;V\;\;\}
 \end{equation}
 In particular then $\g$ is closed under the skew bracket on $\End\,V
 \oplus_{R,T}V$ so that
 $$
  [\;\X\,\oplus_{R,T}0,\;A_x\,\oplus_{R,T}x\;]
  \;\;=\;\;
  \Big(\;[\,\X,\,A_x\,]\;-\;A_{\X\,x}\;\Big)\,\oplus_{R,T}0
  \;+\;A_{\X\,x}\,\oplus_{R,T}\X\,x
 $$
 is necessarily an element of $\g$ proving the congruence $[\X,A_x]\,-\,
 A_{\X x}\,\equiv\,0$ modulo $\h$ for all $\X\,\in\,\h$ and $x\,\in\,V$.
 Similarly the result of the following calculation is an element of $\g$
 \begin{eqnarray*}
  [\;A_x\,\oplus_{R,T}x,\;A_y\,\oplus_{R,T}y\;]
  &=&
  \Big(\;[\,A_x,\,A_y\,]\,-\,R_{x,y}\,-\,A_{A_xy\,-\,A_yx\,-\,T(x,y)}\;\Big)
  \,\oplus_{R,T}0
  \\
  &&
  \qquad+\;A_{A_xy\,-\,A_yx\,-\,T(x,y)}\,\oplus_{R,T}
  \Big(\,A_xy\,-\,A_yx\,-\,T(x,y)\,\Big)
 \end{eqnarray*}
 and thus requires $Q(A,T)_{x,y}\,\equiv\,R_{x,y}$ modulo $\h$ for all
 $x,\,y\,\in\,V$. In passing we remark that $[\,\X\,\oplus_{R,T}0,\,\Y\,
 \oplus_{R,T}0\,]\,=\,[\,\X,\,\Y\,]\,\oplus_{R,T}0$ lies in $\g$ for
 $\X,\,\Y\,\in\,\h$ without further ado.

 \pfill
 In consequence the subspace $\g$ of equation (\ref{galg}) is a skew
 subalgebra of $\End\,V\oplus_{R,T}V$ as soon as $\h\,\star\,A\,\in
 \,V^*\otimes\h$ and $Q(A,T)\,-\,R\,\in\,\L^2V^*\otimes\h$. On the other
 hand the skew bracket on $\End\,V\oplus_{R,T}V$ agrees with the Lie
 bracket of the semidirect product $\End\,V\oplus V$ of $\End\,V$ with
 its representation $V$ up to terms quadratic in $V$. Hence the Jacobiator
 on every skew subalgebra $\g\,\subset\,\End\,V\oplus_{R,T}V$ vanishes
 automatically for two or three arguments in the isotropy subalgebra
 $\h\,\subset\,\End\,V$. In light of this observation the Jacobiator
 satisfies
 \begin{eqnarray*}
  \Jac(\;\X\,\oplus_{R,T}0,\;A_y\,\oplus_{R,T}y,\;A_z\,\oplus_{R,T}z\;)
  &=&
  \Jac(\;\X\,\oplus_{R,T}0,\;0\,\oplus_{R,T}y,\;0\,\oplus_{R,T}z\;)
  \\[2pt]
  &=&
  -\;\Big(\;(\,\X\,\star\,R\,)_{y,\,z}\,\oplus_{R,T}(\,\X\,\star\,T\,)
  (\,y,\,z\,)\;\Big)
 \end{eqnarray*}
 for all $X\,\in\,\h$ and $y,\,z\,\in\,V$, thus the isotropy algebra of
 a Lie subalgebra $\g\,\subset\,\End\,V\oplus_{R,T}V$ projecting onto
 $V$ is necessarily a subalgebra $\h\,\subset\,\stab\,R\,\cap\,\stab\,T$
 of the joint stabilizer of $R$ and $T$ in $\End\,V$. Eventually we
 calculate along a similar line of argument
 \begin{eqnarray*}
  \lefteqn{\Jac(\;A_x\,\oplus_{R,T}x,\,A_y\,\oplus_{R,T}y,\,
  A_z\,\oplus_{R,T}z\;)}\;
  &&
  \\[2pt]
  &=&
  +\;\Jac(\;A_x\,\oplus_{R,T}0,\,0\,\oplus_{R,T}y,\,0\,\oplus_{R,T}z\;)\
  +\;\Jac(\;A_y\,\oplus_{R,T}0,\,0\,\oplus_{R,T}z,\,0\,\oplus_{R,T}x\;)
  \\
  &&
  +\;\Jac(\;A_z\,\oplus_{R,T}0,\,0\,\oplus_{R,T}x,\,0\,\oplus_{R,T}y\;)\;
  +\;\Jac(\;\;\,0\,\;\oplus_{R,T}x,\,0\,\oplus_{R,T}y,\,0\,\oplus_{R,T}z\;)
  \\[2pt]
  &=&
  \Big(\;-\,d^{(\,A,\,T\,)}R\;\Big)_{x,\,y,\,z}\,\oplus_{R,T}
  \Big(\;R\,\wedge\,\id\;-\;d^{(\,A,\,T\,)}T\;\Big)(\,x,\,y,\,z\,)
 \end{eqnarray*}
 using trilinearity, cyclic invariance and the explicit formulas calculated
 above for $\Jac$.
 \qed

 \pfill
 Although the preceeding lemma is reasonably explicit, it is certainly not
 satisfactory in that we would prefer conditions on the triple $(A,R,T)$
 alone, which guarantee the existence of a Lie subalgebra $\g$ of $\End\,V
 \oplus_{R,T}V$ satisfying $\g\,\supset\,\{\;A_x\,\oplus_{R,T}x\;|\;\;x\,
 \in\,V\;\}$. In order to achieve such a reformulation of Lemma \ref{sag}
 let us have another look at the stabilizer filtration
 \begin{equation}\label{df}
  \End\,V\;=\;\ldots\;=\;\h_{-1}\;\supsetneq\;\h_0
  \;\supsetneq\;\ldots\;\supsetneq\;\h_{s-1}
  \;\supsetneq\;\h_s\;=\;\h_{s+1}\;=\;\ldots\;=\;\h_\infty\;\supseteq\;\h
 \end{equation}
 constructed in Definition \ref{sfilt}. Evidently the connection component
 $A\,\in\,V^*\otimes\End\,V$ of the triple $(A,R,T)$ allows us to define
 for every subalgebra $\hat\h\,\subset\,\End\,V$ the derived subalgebra:
 $$
  \hat\h'\;\;:=\;\;
  \{\;\X\,\in\,\hat\h\;\;|\;\;[\;\X,\,A_x\;]\;\equiv\;A_{\X\,x}
  \;\mod\;\hat\h\;\textrm{for all\ }x\,\in\,V\;\;\}
 $$
 For every subalgebra $\hat\h\,\subset\,\End\,V$ the quotient $\End\,V/\hat\h$
 is naturally a representation of $\hat\h$ and the derived subalgebra $\hat\h'$
 is nothing else but the stabilizer of the class $A\,+\,V^*\otimes\hat\h$
 represented by $A$ in $V^*\otimes(\End\,V/\hat\h)$. In particular the
 derived subalgebra is monotone
 $$
  \hat\h_{\mathrm{small}}\;\;\subset\;\;\hat\h_{\mathrm{large}}
  \qquad\Longrightarrow\qquad
  \hat\h'_{\mathrm{small}}\;\;\subset\;\;\hat\h'_{\mathrm{large}}
 $$
 because the canonical projection $V^*\otimes(\End\,V/\hat\h_{\mathrm{small}})
 \longrightarrow V^*\otimes(\End\,V/\hat\h_{\mathrm{large}})$ is
 equivariant under every subalgebra $\hat\h_{\mathrm{small}}\,\subset\,
 \hat\h_{\mathrm{large}}$ so that $\X\,\in\,\hat\h_{\mathrm{small}}$
 stabilizing the class represented by $A$ in $V^*\otimes(\End\,V/\hat
 \h_{\mathrm{small}})$ still stabilizes its image in $V^*\otimes(\End\,V
 /\hat\h_{\mathrm{large}})$.

 \pfill
 Thinking of the derived subalgebra construction as a dynamical
 system $\hat\h\longmapsto\hat\h'$ on the set of subalgebras
 $\hat\h\,\subset\,\End\,V$ we observe that Lemma \ref{sag} is
 actually asking for the isotropy algebra $\h\,=\,\g\,\cap\,\End\,V$
 of a Lie subalgebra $\g\,\subset\,\End\,V\oplus_{R,T}V$ to be a fixed
 point $\h\,=\,\h'$ contained in $\stab\,R\,\cap\,\stab\,T$ while
 at the same time containing all the values of $Q(A,T)\,-\,R\,\in\,
 \L^2V^*\otimes\h$. The latter condition however becomes ever the
 more restrictive the smaller is $\h$, hence the best we can hope
 for is that the unique maximal fixed point $\h_\max$ of the dynamical
 system $\hat\h\longmapsto\hat\h'$ contained in $\stab\,R\,\cap\,\stab\,T$
 is a sufficiently large subalgebra of $\End\,V$ to satisfy $Q(A,T)-R
 \,\in\,\L^2V^*\otimes\h_\max$.

 The existence of this unique maximal fixed point $\h_\max\,\subset\,
 \stab\,R\,\cap\,\stab\,T$ is guaranteed simply by the monotonicity of the
 dynamical system $\hat\h\longmapsto\hat\h'$, to wit $\h_\max\,=\,\h_\infty$
 is necessarily equal to the limit of the filtration sequence (\ref{df})
 starting in $\h_0\,:=\,\stab\,R\,\cap\,\stab\,T$ and iterating $\h_{r+1}
 \,:=\,\h_r'$ the derived subalgebra construction for $r\,\geq\,0$.
 Depending on the curvature--torsion tuple $(R,T)$ this unique maximal
 fixed point $\h_\max$ may or may not satisfy $Q(A,T)-R\,\in\,\L^2V^*
 \otimes\h_\max$. Provided the maximal fixed point $\h_\max$ of the
 dynamical system $\hat\h\longmapsto\hat\h'$ contained in $\stab\,R
 \,\cap\,\stab\,T$ satisfies $Q(A,T)-R\,\in\,\L^2V^*\otimes\h_\max$
 there may of course exist other fixed points $\h\,\subset\,\h_\max$
 still satisfying $Q(A,T)-R\,\in\,\L^2V^*\otimes\h$. In geometric terms
 these additional fixed points correspond to subgroups $G\,\subset\,G_\max$
 still acting transitively on the affine homogeneous space $G_\max/H_\max$.

 In any case the unique maximal fixed point $\h_\max\,=\,\h_\infty$ contained
 in $\stab\,R\,\cap\,\stab\,T$ agrees with the joint stabilizer of all
 iterated covariant derivatives $\nabla^r T$ and $\nabla^r R$ for all
 $r\,\geq\,0$ together according to Lemma \ref{if}. In consequence the
 decisive congruence $Q(A,T)_{x,y}\,\equiv\,R_{x,y}$ modulo $\h_\max$
 is equivalent to the following set of algebraic equations:

 \begin{Theorem}[Algebraic Variety of Affine Homogeneous Spaces]
 \hfill\label{mfh}\break
  For a given connection--curvature--torsion triple $(A,\,R,\,T)$ on a
  vector space $V$ there exists a Lie subalgebra $\g\,\subset\,\End\,V
  \oplus_{R,T}V$ of the skew algebra associated to $(R,T)$ satisfying
  $$
   \g\;\;\supset\;\;
   \im\Big(\;A^{\mathrm{ext}}:\;V\;\longrightarrow\;
   \End\;V\,\oplus_{R,T}V,
   \qquad x\;\longmapsto\;A_x\,\oplus_{R,T}x\;\Big)
  $$
  if and only if $(A,R,T)$ satisfies the first and second Bianchi identities
  of degrees $2$ and $3$
  $$
   d^{(A,T)}T\;\;=\;\;R\,\wedge\,\id
   \qquad\qquad\qquad
   d^{(A,T)}R\;\;=\;\;0
  $$
  as well as the following homogeneous equations of degrees $r+3$ and $r+4$
  for all $r\,\geq\,0$:
  \begin{eqnarray*}
   \Big(\;Q(A,T)\;-\;R\;\Big)\,\bm\,\Big(\;\underbrace{A\,\bm\,(\,A\,\bm\,
   (\,\ldots\,(\,A}_{r\;\,\mathrm{times}}\,\bm\,T\,)\,\ldots\,)\,)\;\Big)
   &=&0
   \\
   \Big(\;Q(A,T)\;-\;R\;\Big)\,\bm\,\Big(\;\underbrace{A\,\bm\,(\,A\,\bm\,
   (\,\ldots\,(\,A}_{r\;\,\mathrm{times}}\,\bm\,R\,)\,\ldots\,)\,)\;\Big)
   &=&0
  \end{eqnarray*}
  The set of all solution triples $(A,R,T)$ to these algebraic equations
  will be denoted $\M(\,\gl\,V\,)$. 
 \end{Theorem}

 \pfill
 Although the algebraic variety $\M(\,\gl\,V\,)$ is formally defined
 by an infinite set of algebraic equations, only a finite number of
 these equations can actually be relevant, after all polynomial rings
 are noetherian rings. Taking the argument leading to these equations
 into account we observe that at least the equations parametrized by
 $r\,\geq\,(\dim\,V)^2$ have to be {\em algebraic} consequences of
 the equations up to $(\dim\,V)^2$. In this observation $(\dim\,V)^2$
 enters simply as a trivial upper bound for the maximal length $s$
 of a sequence $\hat\h_0\,\supsetneq\,\hat\h_1\,\supsetneq\,\ldots\,
 \supsetneq\,\hat\h_s$ of subalgebras of $\End\,V$ with $\hat\h_{r+1}
 \,=\,\hat\h_r'$ for some $A\,\in\,V^*\otimes\End\,V$, which probably
 grows more like $\dim\,V$ than $(\dim\,V)^2$. On the other hand the
 family of examples constructed in Section \ref{fts} shows that there
 is no universal bound independent of $\dim\,V$ for the number of
 algebraic equations needed in order to define $\M(\,\gl\,V\,)$.

 Somewhat more interesting from the theoretical point of view is the
 observation that the algebraic variety $\M(\,\gl\,V\,)$ is actually
 a cone over a projective algebraic variety due to the homogeneity of
 the algebraic equations defining $\M(\,\gl\,V\,)$. More precisely
 the group $\R$ acts via $\lambda\,\star\,(\,A,\,R,\,T\,)\,:=\,
 (\,e^\lambda A,\,e^{2\lambda}R,\,e^\lambda T\,)$ on the set of
 all connection--curvature--torsion triples and thus on the
 invariant subset $\M(\,\gl\,V\,)\,\subset\,V^*\otimes\End\,V\,\times\,
 \L^2V^*\otimes\End\,V\,\times\,\L^2V^*\otimes V$. The same observation
 applies to the algebraic varieties $\M^\tf(\,\gl\,V\,)$ and
 $\M^{\mathrm{ref}}(\,\gl\,V\,)$ of torsion free and reductive
 connection--curvature--torsion triples $(\,A,\,R,\,T\,)$ defined in:

 \begin{Corollary}[Torsion Free Affine Homogeneous Spaces]
 \hfill\label{tf}\break
  Torsion free connection--curvature tuples $(A,R)$ are
  connection--curvature--torsion triples with vanishing torsion $T\,=\,0$.
  Their algebras $\End\,V\oplus_{R,0}V$ allow Lie subalgebras $\g$ with
  $$
   \End\;V\,\oplus_{R,0}V\;\;\supset\;\;\g\;\;\supset\;\;
   \im\Big(\;A^{\mathrm{ext}}:\;\;V\longrightarrow\;\End\,V\,\oplus_{R,0}V,
   \quad x\;\longmapsto\;A_x\oplus_{R,0}x\;\Big)
  $$
  if and only if the curvature $R\,\in\,\L^2V^*\otimes\End\,V$ and the
  connection $A\,\in\,V^*\otimes\End\,V$ satisfy the formally infinite
  system of algebraic equations of degrees $2,\,3$ and $r+4,\,r\,\geq\,0$:
  $$
   R\,\wedge\,\id\;\;=\;\;0
   \qquad
   d^{(A,0)}R\;\;=\;\;0
   \qquad
   \Big(Q(A,0)\,-\,R\Big)\bm\Big(\underbrace{A\bm(A\bm
   (\ldots(A}_{r\;\,\mathrm{times}}\,\bm\,R)\ldots))\Big)
   \;\;=\;\;0
  $$
  The set of all solutions $(A,R)$ to these equations will be denoted
  by $\M^\tf(\,\gl\,V\,)\,\subset\,\M(\,\gl\,V\,)$.
 \end{Corollary}

 \begin{Corollary}[Reductive Affine Homogeneous Spaces]
 \hfill\label{red}\break
  Reductive curvature--torsion tuples $(R,T)$ are
  connection--curvature--torsion triples with vanishing connection
  $A\,=\,0$. The skew algebra $\End\,V\oplus_{R,T}V$ allows Lie subalgebras
  $\g\,\supset\,V$ extending $V$, if and only if the curvature $R\,\in\,
  \L^2V^*\otimes\End\,V$ and torsion $T\,\in\,\L^2V^*\otimes V$ satisfy
  the following system of algebraic equations of degrees $2,\,3,\,3$ and $4$:
  $$
   d^{(0,T)}T\;\;=\;\;R\,\wedge\,\id
   \qquad\quad
   d^{(0,T)}R\;\;=\;\;0
   \qquad\quad
   R\,\bm\,T\;\;=\;\;0
   \qquad\quad
   R\,\bm\,R\;\;=\;\;0
  $$
  In the sequel $\M^{\mathrm{red}}(\,\gl\,V\,)\,\subset\,\M(\,\gl\,V\,)$
  will denote the set of solutions to these equations.
 \end{Corollary}
\section{Additional Parallel Geometric Structures}
\label{par}
 In this short intermediate section we want to discuss variations of the
 algebraic varieties $\M(\,\gl\,V\,)$ parametrizing formal affine homogeneous
 spaces, which take into account additional parallel geometric structures
 like Riemannian metrics and almost complex structures. Because the affine
 geometry determined by the existence of a linear connection of the tangent
 bundle is intrisically a first order geometry, the condition of parallelity
 severely restricts the order of geometric structures we may consider in
 addition, in any case we will restrict ourselves to a discussion of first
 order geometries usually called (first order) $G$--structures. For the
 purpose of this article we prefer to call them $K$--structures instead
 in order to avoid the clash of notation with the big group $G$ of the
 homogeneous space $G/H$.

 \pfill
 In order to define $K$--structures on vector spaces and consequently on
 smooth manifolds $M$ we fix once and for all a closed subgroup $K\,\subset
 \,\GL\,V$ of the general linear group of a vector space $V$ called the model
 space. By definition a $K$--structure on a vector space $T$ of the same
 dimension as $V$ is an equivalence class $\Omega$ of linear isomorphisms
 $F:\,V\longrightarrow T$ under
 $$
  F\;\;\sim_K\;\;\tilde F
  \qquad\Longleftrightarrow\qquad
  F^{-1}\,\circ\,\tilde F\;\;\in\;\;K
 $$
 conversely the linear isomorphisms $F:\,V\longrightarrow T$ representing
 $\Omega$ are called $K$--frames. The general linear group $\GL\,V$ acts
 simply transitively from the right on the set $\Fr(\,V,\,T\,)$ of linear
 isomorphisms $V\longrightarrow T$, thus $K$--structures $\Omega$
 correspond bijectively to points in:
 $$
  \Omega\;\;\in\;\;\Fr(\;V,\;T\;)/_{\displaystyle K}
 $$
 With $K$ being a closed subgroup of $\GL\,V$ by assumption the quotient
 $\Fr(\,V,\,T\,)/K$ is actually a manifold so that we can define a smooth
 $K$--structure on a manifold $M$ of the same dimension as $V$ as a
 smooth section of the quotient bundle
 $$
  \Omega\;\;\in\;\;\Gamma(\;\Fr(\;V,\;TM\;)/_{\displaystyle K}\;)
 $$
 where $\Fr(\,V,\,TM\,)$ is the standard principal $\GL\,V$--bundle of all
 frames on $M$. Such a smooth $K$--structure $\Omega$ is parallel for
 an affine connection $\nabla$ on $TM$, if and only if the parallel
 transport $\mathrm{PT}^\nabla_\gamma:\,T_{\gamma(0)}M\longrightarrow
 T_{\gamma(1)}M$ along arbitrary curves $\gamma:\,[0,1]\longrightarrow M$
 maps $K$--frames to $K$--frames in the sense $\mathrm{PT}^\nabla_{\gamma}
 \circ F\,\in\,\Omega_{\gamma(1)}$ for every $K$--frame $F\,\in\,
 \Omega_{\gamma(0)}$.

 \pfill 
 Whereas the preceeding definition of parallel $K$--structures on a
 manifold $M$ for a closed subgroup $K\,\subset\,\GL\,V$ does not
 comprise the most general geometric structures, nevertheless there
 are quite a number of interesting examples. Consider for example
 $V\,=\,\C^n$ as a real vector space of dimension $2n$, which inherits
 from $\C^n$ the complex structure $I\,\in\,\End\,V$, the real part
 $g\,\in\,\S^2_+V^*$ of the standard hermitean form $(\cdot,\cdot)$
 and the $n$--form $\psi\,\in\,\L^nV^*$:
 $$
  \psi(\;v_1,\,\ldots,\,v_n\;)
  \;\;:=\;\;
  \mathrm{Re}\;\mathrm{det}_\C(\;v_1,\,\ldots,\,v_n\;)
 $$
 The common stabilizer of the triple $(\,g,\,I,\,\psi\,)$ agrees with
 the subgroup $\SU(n)\,\subset\,\GL\,V$ so that an $\SU(n)$--structure
 $\Omega$ on a $2n$--dimensional vector space defines a corresponding triple
 $$
  g_\Omega\;\;:=\;\;g(\,F^{-1}\,\cdot,\,F^{-1}\,\cdot\,)
  \qquad
  I_\Omega\;\;:=\;\;F\,\circ\,I\,\circ\,F^{-1}
  \qquad
  \psi_\Omega\;\;:=\;\;\psi(\,F^{-1}\,\cdot,\,\ldots,\,F^{-1}\,\cdot\,)
 $$
 on $T$ for any representative $F\,\in\,\Omega$. On the other hand the
 model triple $(\,g,\,I,\,\psi\,)$ satisfies
 \begin{equation}\label{su}
  I^2\;\;=\;\;-\,\id_V
  \qquad
  g(\,I\,\cdot,\,I\,\cdot\,)\;\;=\;\;g
  \qquad
  \mathrm{Der}^2_I\,\psi\;\;=\;\;-\,n^2\,\psi
  \qquad
  g^{-1}(\,\psi,\,\psi\,)\;\;=\;\;2^{n-1}
 \end{equation}
 and it can be shown straightforwardly that $\GL\,V$ acts transitively on
 the set of solutions to these algebraic equations in $\S^2_+V^*\,\times\,
 \End\,V\,\times\,\L^nV^*$. In consequence $\SU(n)$--structures $\Omega$
 on a vector space $T$ are in bijection to triples $(\,g,\,I,\,\psi\,)$
 in $\S^2_+T^*\,\times\,\End\,T\,\times\,\L^nT^*$ satisfying the algebraic
 equations (\ref{su}). Similarly the tuple $(\,g,\,I\,)$ defines the
 underlying $\mathbf{U}(n)$--structure and $g\,\in\,\S^2_+T^*$ only
 the underlying $\mathbf{O}(2n)$--structure on $T$.

 \pfill
 Coming back to homogeneous spaces $G/H$ we define a left invariant
 $K$--structure as a left invariant section of the homogeneous quotient
 bundle $\Fr(\,V,\,T(G/H)\,)/K$ of the homogeneous frame bundle. By the
 classification of left invariant sections such a section $\Omega$ is
 completely determined by its $H$--invariant value in the base point
 $eH\,\in\,G/H$ of $G/H$
 $$
  \Omega_{eH}\;\;\in\;\;\Big[\;\Fr(\;V,\;\g/\h\;)/_{\displaystyle K}\;\Big]^H
 $$
 where the condition of $H$--invariance reads for a representative
 $F:\,V\longrightarrow\g/\h$ of $\Omega_{eH}$
 $$
  (\,h\,\star\,)\,\circ\,F\;\;\sim_K\;\;F
  \qquad\Longleftrightarrow\qquad
  F^{-1}\,\circ\,(\,h\,\star\,)\,\circ\,F\;\;\in\;\;K
 $$
 for every $h\,\in\,H$, equivalently the image $\mathrm{Ad}\,H
 \,\subset\,\GL\,\g/\h$ of $H$ under the adjoint representation must
 be conjugated under $F$ to the subgroup $F^{-1}(\,\mathrm{Ad}\,H\,)F
 \,\subset\,K$ of $K$. Similarly a left invariant $K$--structure on
 $G/H$ is parallel for a left invariant connection $A:\,\g\longrightarrow
 \End\,\g/\h$ on $T(G/H)$ provided $F^{-1}\,\circ\,A_X\,\circ\,F$ is
 an element of the Lie algebra $\k$ of $K$ for every $X\,\in\,\g$.

 In consequence we may define a formal affine homogeneous space
 with parallel $K$--structure as a pair of Lie algebras $\g\,
 \supset\,\h$ endowed with a formal left in variant connection
 $A:\,\g\longrightarrow\End\,\g/\h$ and a $K$--structure $\Omega$
 on $\g/\h$ such that $F^{-1}\circ(\,H\,\star\,)\circ F\,\subset\,\k$
 and $F^{-1}\circ A_X\circ F\,\in\,\k$ for all $H\,\in\,\h$ and
 $X\,\in\,\g$. On the other hand we recall that we need to choose
 a frame $F:\,V\longrightarrow\g/\h$ anyhow in order to associate
 to such a formal homogeneous space a point in the algebraic variety
 $\M(\,\gl\,V\,)$, making the straightforward choice of a frame
 $F$ representing $\Omega$ we arrive at the algebraic variety of
 formal affine homogeneous spaces
 \begin{equation}\label{msk}
  \M(\;\k\;)
  \;\;:=\;\;
  \M(\;\gl\;V\;)\;\cap\;(\;V^*\,\otimes\,\k\;\times\;
  \L^2V^*\,\otimes\,\k\;\times\;\L^2V^*\,\otimes\,V\;)
 \end{equation}
 with parallel $K$--structure. Passing through the arguments used
 in the construction of the algebraic variety $\M(\,\gl\,V\,)$ we
 conclude that the connection--curvature--torsion triples $(A,R,T)$
 in $\M(\;\k\;)$ correspond to maximal Lie subalgebras $\g\,\subset
 \,\k\,\oplus_{R,T}V$ of the skew algebras $\k\,\oplus_{R,T}V$, which
 contain the image of the extended connection $V\longrightarrow\k\,
 \oplus_{R,T}V,\,x\longmapsto A_x\,\oplus_{R,T}x$. Simi\-lar considerations
 apply of course to the algebraic varieties $\M^\tf(\,\k\,)$ and
 $\M^{\mathrm{red}}(\,\k\,)$ of torsion free and reductive formal
 affine homogeneous spaces with parallel $K$--structures.
\section{Contact Order and Formal Tangent Space}
\label{fts}
 In a sense the algebraic variety $\M(\,\gl\,V\,)$ of affine homogeneous spaces
 and its variants are only coarse moduli spaces, because they parametrize
 isometry classes of formal affine homogeneous spaces $\g\,\supset\,\h$
 augmented by a frame isomorphism $F:\,V\longrightarrow\g/\h$ and a split
 $\g/\h\longrightarrow\g$ of the canonical projection. The effect of changing
 the frame and/or the split introduces an equivalence relation on $\M(\,\gl
 \,V\,)$, the equivalence relation $\sim_\infty$ of contact to all orders,
 which is approximated by the equivalence relations $\sim_d$ of contact
 to order $d\,\in\,\N_0$, In turn the formal tangent space to the true moduli
 space $\M_\infty(\,\gl\,V\,)\,:=\,\M(\,\gl\,V\,)/\sim_\infty$ in a point
 $[A,\,R,\,T]$ is filtered by directions staying in contact with
 $[A,\,R,\,T]$ up to order $d\,\geq\,0$.

 After introducing the contact equivalence relations $\sim_d$ and
 $\sim_\infty$ we use the stabilizer filtration of Definition \ref{sfilt}
 in order to associate to every connection--curvature--torsion triple
 $(A,\,R,\,T)\,\in\,\M(\,\gl\,V\,)$ its Spencer cohomology $H^{\bullet,\circ}
 (\,\h\,)$. Moreover we identify this Spencer cohomology with the sucessive
 filtration quotients of the contact filtration on the formal tangent space
 $T_{[A,R,T]}\M_\infty(\,\gl\,V\,)$ to the true moduli space $\M_\infty
 (\,\gl\,V\,)$ in the point $[A,\,R,\,T]$. Eventually we will illustrate
 this interpretation of the Spencer cohomology associated to $(A,\,R,\,T)$
 by detailed calculations for the family of pairwise non--isometric
 Riemannian homogeneous spaces with large Singer invariant constructed
 by Meusers \cite{m}.

 \begin{Definition}[Contact Relation for Connection--Curvature--Torsion
  Triples]
 \hfill\label{jet}\break
  Two connection--curvature--torsion triples $(A,R,T)$ and $(\tilde A,
  \tilde R,\tilde T)$ on a finite--dimensional vector space $V$ are said
  to be in contact to order $d\,\geq\,0$ written $(A,R,T)\,\sim_d\, 
  (\tilde A,\tilde R,\tilde T)$, if there exists a linear automorphism
  $F\,\in\,\GL\,V$ of $V$ pulling the formal covariant derivatives of
  $\tilde R$ and $\tilde T$ of all orders $r\,=\,0,\,\ldots,\,d$ back
  to the formal covariant derivatives of $R$ and $T$:
  \begin{eqnarray*}
   F^*\Big(\,\underbrace{\tilde A\,\bm\,(\,\tilde A\,\bm\,\ldots\,
   (\,\tilde A}_{r\;\,\mathrm{times}}\,\bm\,\tilde T\,)\,\ldots\,)\,\Big)
   &=&
   \underbrace{A\,\bm\,(\,A\,\bm\,\ldots\,(\,A}_{r\;\,\mathrm{times}}
   \,\bm\,T\,)\,\ldots\,)
   \\
   F^*\Big(\,\underbrace{\tilde A\,\bm\,(\,\tilde A\,\bm\,\ldots\,
   (\,\tilde A}_{r\;\,\mathrm{times}}\,\bm\,\tilde R\,)\,\ldots\,)\,\Big)
   &=&
   \underbrace{A\,\bm\,(\,A\,\bm\,\ldots\,(\,A}_{r\;\,\mathrm{times}}
   \,\bm\,R\,)\,\ldots\,)
  \end{eqnarray*}
  Similar the notation $(A,R,T)\,\sim_\infty\,(\tilde A,\tilde R,\tilde T)$
  indicates triples in contact to all orders $d\,\geq\,0$.
 \end{Definition}

 \pfill
 The naturality of the operation $\bm$ implies of course for the iterated
 covariant derivatives
 $$
  F^*(\,A\,\bm\,(\,A\,\bm\,\ldots\,(\,A\,\bm\,T\,)\,\ldots\,)\,\Big)
  \;\;=\;\; 
  F^*A\,\bm\,(\,F^*A\,\bm\,\ldots\,(\,F^*A\,\bm\,F^*T\,)\,\ldots\,)
 $$
 of $T$ or similarly of $R$. The main argument of Lemma \ref{if}
 may thus be varied to prove:

 \begin{Lemma}[Explicit Form of the Contact Relation]
 \hfill\label{cjet}\break
  Two connection--curvature--torsion triples $(A,R,T)$ and $(\tilde A,
  \tilde R,\tilde T)$ on a finite--dimensional vector space $V$ are in
  contact $(A,R,T)\,\sim_d\,(\tilde A,\tilde R,\tilde T)$ to order
  $d\,\geq\,0$, if and only if there exists a linear automorphism
  $F\,\in\,\GL\,V$ satisfying $F^*\tilde T\,=\,T$ and $F^*\tilde R\,=\,R$
  as well as:
  $$
   F^*\tilde A
   \;\;\equiv\;\;
   A \qquad\mod\;\;\;V^*\,\otimes\,\h_{d-1}
  $$
  In consequence two triples $(A,R,T)\,\sim_d\,(\tilde A,\tilde R,\tilde T)$
  in contact to an order $d\,>\,\Singer(A,R,T)$ exceeding the Singer invariant
  of one are in contact $(A,R,T)\,\sim_\infty\,(\tilde A,\tilde R,\tilde T)$
  to all orders.
 \end{Lemma}

 \pfill
 The infinite order contact relation $\sim_\infty$ reflects exactly the
 dependence of the connection--curvature--torsion triple $(A,R,T)$ associated
 to a formal affine homogeneous space $\g\,\supset\,\h$ on the additional
 choice of frame $F:\,V\longrightarrow\g/\h$ and split $\g/\h\longrightarrow
 \g$ of the canonical projection. The moduli space of isometry classes of
 formal affine homogeneous spaces may thus be defined as the set of
 equivalence classes of points $(A,R,T)\,\in\,\M(\,\gl\,V\,)$ under
 $\sim_\infty$:
 $$
  \M_\infty(\;\gl\,V\;)\;\;:=\;\;\M(\;\gl\,V\;)/_{\displaystyle\sim_\infty}
 $$
 Similar definitions can be made for the moduli spaces of torsion free
 $\M^\tf_\infty(\,\gl\,V\,)$ or reductive formal affine homogeneous spaces
 $\M^\mathrm{red}_\infty(\,\gl\,V\,)$ with or without additional left invariant
 parallel $K$--structures. It turns out that the deformation theory of
 isometry classes of formal affine homogeneous spaces considered as points
 in $\M_\infty(\,\gl\,V\,)$ is governed by a suitable version of the Spencer
 cohomology associated to the purely algebraic concept of a graded comodule
 over the symmetric coalgebra $\S\,V^*$ of a vector space $V$:
 
 \begin{Definition}[Comodules over Symmetric Coalgebras]
 \hfill\label{comod}\break
  A comodul over the symmetric coalgebra $\S\,V^*$ of a finite--dimensional
  vector space $V$ is a $\Z$--graded vector space $\A^\bullet$ together
  with a bilinear map called directional derivative $V\times\A^\bullet
  \longrightarrow\A^{\bullet-1},\,(y,\X)\longmapsto\frac{\partial\X}
  {\partial y},$ homogeneous of degree $-1$ with respect to $\A$ such that
  $$
   \frac{\partial}{\partial y}\Big(\;\frac{\partial\X}{\partial z}\;\Big)
   \;\;=\;\;
   \frac{\partial}{\partial z}\Big(\;\frac{\partial\X}{\partial y}\;\Big)
  $$
  for all $y,\,z\,\in\,V$. In particular the iterated directional derivatives
  like $\frac{\partial^2\X}{\partial y\,\partial z}$ are well defined.
 \end{Definition}

 \begin{Definition}[Spencer Complex and Cohomology]
 \hfill\label{scc}\break
  The Spencer complex associated to a comodule $\A$ over the symmetric
  coalgebra $\S\,V^*$ of  a finite--dimensional vector space $V$ is the
  bigraded complex of alternating multilinear forms $\L^\circ V^*\otimes
  \A^\bullet$ on $V$ with values in the comodule $\A$ 
  $$
   \ldots\;\stackrel B\longrightarrow\;
   \L^{\circ-1}V^*\,\otimes\,\A^{\bullet+1}
   \stackrel B\longrightarrow\;
   \L^\circ V^*\,\otimes\,\A^\bullet
   \stackrel B\longrightarrow\;
   \L^{\circ+1}V^*\,\otimes\,\A^{\bullet-1}
   \stackrel B\longrightarrow\;\ldots
  $$
  endowed with the Spencer coboundary operator $B$ defined for
  $\eta\,\in\,\L^kV^*\otimes\A^\bullet$ by:
  $$
   (\,B\,\eta\,)(\;x_0,\,\ldots,\,x_k\;)
   \;\;:=\;\;
   \sum_{\mu\,=\,0}^k(\,-1\,)^\mu\;
   \frac{\partial}{\partial x_\mu}
   B(\;x_0,\,\ldots,\,\widehat{x_\mu},\,\ldots,\,x_k\;)
  $$
  The corresponding bigraded cohomology theory for comodules is called
  Spencer cohomology:
  $$
   H^{\bullet,\circ}(\;\A\;)
   \;\;:=\;\;
   \frac{\ker(\;B:\;
   \hbox to28pt{\hfill$\A^\bullet$\hfill}\otimes
   \hbox to38pt{\hfill$\L^\circ V^*$\hfill}\;\longrightarrow\;
   \hbox to28pt{\hfill$\A^{\bullet-1}$\hfill}\otimes
   \hbox to38pt{\hfill$\L^{\circ+1}V^*$\hfill}\;)}
   {\im(\;\,B:\,\;
   \hbox to28pt{\hfill$\A^{\bullet+1}$\hfill}\otimes
   \hbox to38pt{\hfill$\L^{\circ-1}V^*$\hfill}\;\longrightarrow\;
   \hbox to28pt{\hfill$\A^\bullet$\hfill}\otimes
   \hbox to38pt{\hfill$\L^\circ V^*$\hfill}\;)}
  $$
 \end{Definition}

 \pfill
 A detailed description of the general properties of the Spencer cohomology
 of comodules is certainly beyond the scope of this article, for more
 information see for example \cite{bcg},\ \cite{exw}. Nevertheless we want
 to point out that the Spencer cohomology $H^{\bullet,\circ}(\,\A\,)$ of a
 comodule $\A^\bullet$ is naturally a graded right module over the exterior
 algebra $\L^\circ V^*$. For a comodule $\A$ constant in the directions of
 a subspace $W\,\subset\,V$ in the sense $\frac{\partial\X}{\partial x}\,=\,0$
 for all $\X\,\in\,\A$ and all $x\,\in\,W$ for example the right
 $\L^\circ V^*$--module structure turns out to be convenient to prove
 \begin{equation}\label{sc}
  H^{\bullet,\,\circ}(\;\A\;)
  \;\;\cong\;\;
  H^{\bullet,\,\circ}_{V/W}(\;\A\;)
  \,\otimes\,\L^\circ W^*
 \end{equation}
 where $H^{\bullet,\circ}_{V/W}(\,\A\,)$ refers to the Spencer cohomology
 of $\A^\bullet$ considered as a comodule over the coalgebra $\S\,(V/W)^*$
 and the repeated grading symbol $\circ$ indicates the product grading.

 \begin{Definition}[Formal Directional Derivatives]
 \hfill\label{fdd}\break
  Consider the strictly decreasing filtration (\ref{df}) associated to
  a connection--curvature--torsion triple $(A,R,T)$ in the variety
  $\M(\,\gl\,V\,)$. The direct sum of sucessive filtration quotients
  $$
   \h^\bullet
   \;\;=\;\;
   \bigoplus_{r\,\in\,\Z}
   \Big(\;\h_{r-1}/_{\displaystyle\h_r}\;\Big)
   \;\;=\;\;
   (\,\End\,V/_{\displaystyle\h_0}\,)
   \;\oplus\;\ldots\;\oplus\;
   (\,\h_{s-1}/_{\displaystyle\h_\infty}\,)
  $$
  is a comodule over the symmetric coalgebra $\S\,V^*$ under directional
  derivatives defined by
  $$
   \frac{\partial\X}{\partial x}
   \;\;:=\;\;
   [\;\X,\,A_x\;]\;-\;A_{\X\,x}\;+\;\h_{r-1}
  $$
  for $x\,\in\,V$ and all representatives $\X\,\in\,\h_{r-1}$ of a class
  $\X\,+\,\h_r$ in the quotient $\h^r\,:=\,\h_{r-1}/\h_r$.
 \end{Definition}

 \pfill
 Of course we should not forget to verify the axiomatic commutation
 of the formal directional derivatives for a comodule over the symmetric
 coalgebra $\S\,V^*$, which is quite surprising in view of the complicated
 definition of the comodule $\h^\bullet$ and its formal directional
 derivatives. Straightforward calculation of iterated formal derivatives
 results in the not too pleasant
 \begin{eqnarray*}
  \frac\partial{\partial y}\,\frac{\partial\X}{\partial z}
  &\equiv&
  [\;[\,\X,\,A_z\,]\,-\,A_{\X z},\;A_y\;]
  \;-\;A_{(\,[\,\X,\,A_z\,]\;-\;A_{\X z}\,)\,y}
  \\
  &\equiv&
  -\;[\;A_y,\;[\,\X,\,A_z\,]\;]\;+\;[\;A_y,\;A_{\X z}\;]
  \;-\;A_{\X(\,A_zy\,)}\;+\;A_{A_z(\,\X y\,)}\;+\;A_{A_{\X z}y}
 \end{eqnarray*}
 modulo $\h_{r-2}$ for a given representative $\X\,\in\,\h_{r-1}$ of
 a class $\X\,+\,\h_r\,\in\,\h_{r-1}/\h_r$, in consequence
 \begin{eqnarray*}
  \frac\partial{\partial z}\,\frac{\partial\X}{\partial y}
  \;-\;
  \frac\partial{\partial y}\,\frac{\partial\X}{\partial z}
  &\equiv&
  +\;[\;\X,\;[\,A_y,\,A_z\,]\;-\;A_{A_yz}\;+\;A_{A_zy}\;+\;A_{T(y,z)}\;]
  \;-\;A_{\X\,T(y,z)}
  \\[-2pt]
  &&
  -\;[\;A_{\X y},\;A_z\;]\;+\;A_{A_{\X y}z}\;-\;A_{A_z(\,\X y\,)}
  \;-\;A_{T(\X y,z)}\;+\;A_{T(\X y,z)}
  \\[2pt]
  &&
  -\;[\;A_y,\;A_{\X z}\;]\;+\;A_{A_y(\,\X z\,)}\;-\;A_{A_{\X z}y}
  \;-\;A_{T(y,\X z)}\;+\;A_{T(y,\X z)}
  \\[2pt]
  &&
  +\;[\;\X,\;A_{A_yz}\;-\;A_{A_zy}\;-\;A_{T(y,z)}\;]
  \;-\;A_{\X(\,A_yz\;-\;A_zy\;-\;T(y,z)\,)}
  \\
  &\equiv&
  \Big[\;\X\,\star\,Q(\,A,\,T\,)\;\Big]_{y.\,z}
  \;+\;
  \frac{\partial\X}{\partial(\,A_yz\,-\,A_zy\,-\,T(y,\,z)\,)}
  \;-\;
  A_{[\,\X\,\star\,T\,](\,y,\,z\,)}
 \end{eqnarray*}
 modulo $\h_{r-2}$. For $r\,\leq\,1$ there is actually nothing to prove,
 because the right hand side vanishes for trivial reasons due to $\h_{-1}
 \,=\,\End\,V$. On the other hand $\X\,\in\,\h_{r-1}\,\subset\,\h_0$ for
 $r\,\geq\,2$ so that $\X\,\star\,T\,=\,0$ by the very definition of
 $\h_0\,=\,\stab\,R\,\cap\,\stab\,T$. The directional derivative term
 on the right hand side vanishes by definition modulo $\h_{r-2}$ leaving
 us with
 $$
  \frac\partial{\partial z}\,\frac{\partial\X}{\partial y}
  \;-\;
  \frac\partial{\partial y}\,\frac{\partial\X}{\partial z}
  \;\;\equiv\;\;
  \Big[\;\X\,\star\,(\;Q(\,A,\,T\,)\;-\;R\;)\;\Big]_{y,\,z}
  \qquad
  \mod\;\;\h_{r-2}
 $$
 in view of $\X\,\in\,\stab\,R$. For a connection--curvature--torsion triple
 $(A,R,T)\,\in\,\M(\,\gl\,V\,)$ however the difference $Q(A,T)\,-\,R$ is an
 $\h_\infty$--valued $2$--form on $V$ so that $\X\,\star\,(\,Q(A,T)\,-\,R\,)$
 is actually $\h_{r-1}$--valued for all $\X\,\in\,\h_{r-1}\,\supset\,
 \h_\infty$. In consequence the formal directional derivatives of the
 comodule $\h^\bullet$ defined above do in fact commute and $\h^\bullet$
 is a comodule over the symmetric coalgebra $\S\,V^*$. Its Spencer cohomology
 in form degree $0$ is very simple:

 \begin{Remark}[Stabilizer Filtration and Spencer Cohomology]
 \hfill\label{sfsc}\break
  The only non--vanishing Spencer cohomology in form degree $0$ of the
  comodule $\h^\bullet$ associated to a connection--curvature--torsion
  triple $(A,R,T)$ in the algebraic variety $\M(\,\gl\,V\,)$ is:
  $$
   H^{0,0}(\;\h\;)\;\;=\;\;\h^0\;\;:=\;\;\End\,V/_{\displaystyle\h_0}
  $$
  In fact $H^{r,0}(\,\h\,)\,=\,\{0\}$ for all positive $r\,>\,0$ by the
  definition the stabilizer filtration (\ref{df}).
 \end{Remark}

 \pfill
 Whereas the Spencer cohomology in form degree $0$ is thus not particularly
 interesting, the Spencer cohomology $H^{\bullet,1}(\,\h\,)$ in form degree
 $1$ has the following neat geometric interpretation in terms of deformations
 of a formal affine homogeneous space $\g\,\supset\,\h$.

 \begin{Lemma}[Formal Tangent Space of the Moduli Space]
 \hfill\label{modts}\break
  Given a point $[\,A,\,R,\,T\,]$ in the moduli space $\M_\infty(\,V\,)$
  of formal affine homogeneous spaces of dimension $\dim\,V$ we may
  consider the subsets of formal affine homogeneous spaces
  $$
   \M_\infty^d(\;A,\,R,\,T\;)
   \;\;:=\;\;
   \{\;[\,\tilde A,\,\tilde R,\,\tilde T\,]\;\;|\;\;
   [\,\tilde A,\,\tilde R,\,\tilde T\,]\,\sim_d\,[\,A,\,R,\,T\,]
   \;\}
   \;\;\subset\;\;
   \M_\infty(\,\gl\,V\,)
  $$
  in contact to $[\,A,\,R,\,T\,]$ to order $d\,\geq\,0$. The resulting
  decreasing filtration of $\M_\infty(\,\gl\,V\,)$
  $$
   \M_\infty(\,\gl\,V\,)
   \;\supseteq\;
   \M_\infty^0(A,R,T)
   \;\supseteq\;
   \M_\infty^1(A,R,T)
   \;\supseteq\;\ldots\;\supseteq\;
   \M_\infty^s(A,R,T)
   \;\supseteq\;
   \{\,[A,R,T]\,\}
  $$
  provides an interpretation of the Spencer cohomology of $\h$ by means
  of formal tangent spaces:
  $$
   H^{\bullet,\,1}(\;\h\;)
   \;\;=\;\;
   T_{[A,R,T]}\M^\bullet_\infty(A,R,T)
   /_{\displaystyle T_{[A,R,T]}\M^{\bullet-1}_\infty(A,R,T)}
  $$
 \end{Lemma}

 \proof
 Consider a curve of connection--curvature--torsion triples
 $\varepsilon\longmapsto(\,A_\varepsilon,\,R_\varepsilon,\,T_\varepsilon\,)$
 representing a tangent vector to $\M_\infty^d(A,R,T)$ in the point
 $[\,A,\,R,\,T\,]\,=\,[\,A_0,\,R_0,\,T_0\,]$:
 $$
  \left.\frac d{d\varepsilon}\right|_0
  [\;A_\varepsilon,\;R_\varepsilon,\;T_\varepsilon\;]
  \;\;\in\;\;T_{[\,A,\,R,\,T\,]}\M_\infty^d(A,R,T)
 $$
 By assumption all triples $[A_\varepsilon,\,R_\varepsilon,\,T_\varepsilon]$
 are in contact to order $d\,\geq\,0$ so that we may find a smooth
 curve $\varepsilon\longmapsto F_\varepsilon$ in $\GL\,V$ satisfying
 $F^*_\varepsilon T_\varepsilon\,=\,T$ and $F^*_\varepsilon R_\varepsilon
 \,=\,R$ as well as $F^*_\varepsilon A_\varepsilon\,\equiv\,A$ modulo
 $V^*\otimes\h_{d-1}$. The infinitesimal variation of the curve $\varepsilon
 \longmapsto[\,A_\varepsilon,\,R_\varepsilon,\,T_\varepsilon\,]$ defined
 by
 $$
  \delta\,A
  \;\;:=\;\;
  \left.\frac d{d\varepsilon}\right|_0
  F^*_\varepsilon A_\varepsilon
  \;\;\in\;\;
  V^*\,\otimes\,\h_{d-1}
 $$
 clearly satisfies $\delta A\,\in\,V^*\otimes\h_d$, if the original curve
 $\varepsilon\longmapsto(\,A_\varepsilon,\,R_\varepsilon,\,T_\varepsilon\,)$
 stays actually in contact to $[\,A,\,R,\,T\,]$ to order $d+1$. For this
 reason we are interested in the class represented by $\delta A$ in
 $V^*\otimes\h^d\,=\,V^*\otimes(\h_{d-1}/\h_d)$. This infinitesimal
 variation class is closed, because
 \begin{eqnarray*}
  B(\,\delta A\,)(\,y,z\,)
  &\equiv&
  \frac\partial{\partial y}(\,\delta A_z\,)
  \,-\,\frac\partial{\partial z}(\,\delta A_y\,)
  \;\;\equiv\;\;
  [\,\delta A_z,\,A_y\,]\,-\,[\,\delta A_y,\,A_z\,]
  \,+\,A_{\delta A_yz\,-\,\delta A_zy}
  \\[2pt]
  &\equiv&
  -\;\delta Q(\,A,\,T\,)_{y,\,z}\;-\;\delta A_{A_yz\,-\,A_zy}
 \end{eqnarray*}
 vanishes modulo $\h_{d-1}$, where $\delta Q$ is the infinitesimal
 variation of the approximate curvature:
 \begin{equation}\label{mainc}
  \delta Q\;\;:=\;\;
  \left.\frac d{d\varepsilon}\right|_0
  Q(\;F^*_\varepsilon A_\varepsilon,\;T\;)
  \;\;\in\;\;\L^2V^*\otimes\h_{d-1}
 \end{equation}
 In fact all connection--curvature--torsion triples $(\,F^*_\varepsilon
 A_\varepsilon,\,R,\,T\,)\,\in\,\M(\,\gl\,V\,)$ are in contact to order
 $d\,\geq\,0$ by assumption and thus share the initial part $\h_{-1}\,
 \supsetneq\,\ldots\,\supsetneq\,\h_{d-1}$ of their stabilizer filtrations
 by Lemma \ref{if}. Among the equations defining $\M(\,\gl\,V\,)$ is the
 congruence
 $$
  Q(\;F^*_\varepsilon A_\varepsilon,\;T\;)\;\equiv\;R
  \;\;\mod\;\;\L^2V^*\,\otimes\,(\,\h_\varepsilon\,)_\infty
 $$
 which implies $Q(F^*_\varepsilon A_\varepsilon,T)\,\equiv\,R$ modulo
 $\L^2V^*\otimes\h_{d-1}$ for all $\varepsilon$ and thus proves the
 implicit claim in definition (\ref{mainc}). Being closed for the Spencer
 coboundary operator $B$ the infinitesimal variation represents a class
 $[\,\delta A\,]\,\in\,H^{d,1}(\,\h\,)$ in the Spencer cohomology of the
 comodule $\h^\bullet$ associated to $[\,A,\,R,\,T\,]$. In case this class
 vanishes and the infinitesimal variation is exact
 $$
  \delta\,A\;\;=\;\;B\,\X\;\;\equiv\;\;
  \X\,\star\,A\;\;\mod\;\;V^*\,\otimes\,\h_d
 $$
 for some $\X\,\in\,\h_d$, hence the two curves $\varepsilon\longmapsto
 (\,e^{-\varepsilon\X}\,\star\,F^*_\varepsilon A_\varepsilon,\,R,\,T\,)$ and
 $\varepsilon\longmapsto(\,A_\varepsilon,\,R_\varepsilon,\,T_\varepsilon\,)$
 represent the same tangent vector in $[\,A,\,R,\,T\,]$ tangent to
 $\M^{d+1}_\infty(\,\gl\,V\,)$.
 \qed

 \pfill
 In order to illustrate the relation between the Spencer cohomology
 of Definition \ref{scc} and the formal tangent spaces to the moduli
 space $\M_\infty(\,\gl\,V\,)$ we want to study the family of examples
 of Riemannian homogeneous spaces with large Singer invariant constructed
 by Meusers \cite{m} in detail. Calculations can be streamlined significantly
 using the standard identification of the Lie algebra $\so\,V$ of skew
 symmetric endomorphisms on a euclidian vector space $V$ with scalar
 product $g$ with the second exterior power $\L^2V$
 $$
  \Lambda^2V\;\stackrel\cong\longrightarrow\;\so\;V,
  \qquad X\,\wedge\,Y\;\longmapsto\;
  \Big(\;Z\;\longmapsto\;g(\,X,\,Z\,)\,Y\;-\;g(\,Y,\,Z\,)\,X\;\Big)
 $$
 characterized by $g(\,\X,\,Y\wedge Z\,)\,=\,g(\,\X Y,\,Z\,)$ for all
 $Y,\,Z\,\in\,V$ and $\X\,\in\,\L^2V$ on the left, but $\X\,\in\,\so\,V$
 on the right hand side. The Lie bracket on $\L^2V\,=\,\so\,V$ satisfies
 the formulas
 $$
  [\;X\wedge Y,\;X\wedge Z\;]
  \;\;=\;\;
  g(X,X)\;Y\wedge Z
  \qquad\qquad
  [\;\X\,,\;Y\wedge Z\;]
  \;\;=\;\;
  \X\,Y\wedge Z\;+\;Y\wedge\X\,Z
 $$
 for all $\X\,\in\,\so\,V$ and all $Y,\,Z\,\in\,V$ satisfying $g(X,Y)\,=\,
 0\,=\,g(X,Z)$. Consider now a euclidian vector space $V_\circ$ endowed
 with a scalar product $g\,\in\,\S^2_+V^*$ and an endomorphism $F:\,
 V_\circ\longrightarrow V_\circ$. The direct sum $V\,:=\,\R\oplus V_\circ$
 is then a euclidian Lie algebra under the extension $g(\,x\oplus X,\,
 y\oplus Y\,)\,=\,xy\,+\,g(X,Y)$ of the scalar product and the Lie bracket:
 $$
  [\;x\,\oplus\,X,\;y\,\oplus\,Y\;]
  \;\;:=\;\;
  0\,\oplus\,(\;x\,FY\;-\;y\,FX\;)
 $$
 Evidently $V$ is a solvable Lie algebra with abelian nilpotent subalgebra
 $[\,V,\,V\,]\,=\,V_\circ$, which can be realized as the Lie algebra of left
 invariant vector fields on the corresponding simply connected solvable
 Lie group $G$. After a short calculation it turns out that the Levi--Civita
 connection of the left invariant metric $g$ on $G$ depends on the
 decomposition of $F\,=\,F_+\,+\,F_-$ into its symmetric part $F_+$
 and its skew symmetric part $F_-$, more precisely we obtain:
 $$
  A:\quad V\;\longrightarrow\;\End\,V,\qquad
  x\,\oplus\,X\;\longrightarrow\;(\,F_+X\,)\,\wedge\,\1\;+\;x\,F_-
 $$
 In fact $A_{x\,\oplus\,X}\,\in\,\so\,V$ is a skew symmetric endomorphism
 of $V$ for all $x\oplus X\,\in\,V$ and
 \begin{eqnarray*}
  \lefteqn{A_{x\,\oplus X}(\,y\,\oplus\,Y\,)
   \;-\;A_{y\,\oplus Y}(\,x\,\oplus\,X\,)}
  &&
  \\[2pt]
  &=&
  \Big(\,g(\,F_+X,\,Y\,)\,-\,g(\,F_+Y,\,X\,)\,\Big)\oplus
  \Big(\,(\,-\,y\,F_+X\,+\,x\,F_-Y\,)\,-\,(\,-\,x\,F_+Y\,+\,y\,F_-X\,)\,\Big)
  \\
  &=&
  0\,\oplus\,\Big(\;x\,FY\;-\;y\,FX\;\Big)
  \;\;=\;\;
  [\;x\,\oplus\,X,\;y\,\oplus\,Y\;]
 \end{eqnarray*}
 agrees with the Lie bracket. Using this piece of information we calculate
 the curvature to be
 \begin{eqnarray*} 
  \lefteqn{R_{x\,\oplus\,X,\,y\,\oplus\,Y}}
  &&
  \\[4pt]
  &=&
  [\;(\,F_+X\,)\,\wedge\,\1\;+\;x\,F_-,\;
  (\,F_+Y\,)\,\wedge\,\1\;+\;y\,F_-\;]
  \;-\;
  (\;x\,F_+F\,Y\;-\;y\,F_+F\,X\;)\,\wedge\,\1
  \\[4pt]
  &=&
  (\,F_+X\,)\,\wedge\,(\,F_+Y\,)
  \;-\;x\,(\,F_+^2\;+\;[\,F_+,\,F_-\,]\;)\,Y\,\wedge\,\1
  \;+\;y\,(\,F_+^2\;+\;[\,F_+,\,F_-\,]\;)\,X\,\wedge\,\1
 \end{eqnarray*}
 In order to determine the stabilizer of $R$ in $\so\,V$ it seems prudent
 to study the Ricci curvature:
 $$
  \mathrm{Ric}(\;y\,\oplus\,Y,\;z\,\oplus\,Z\;)
  \;\;:=\;\;
  \tr_V\Big(\;x\,\oplus\,X\;\longmapsto\;
  R_{x\,\oplus\,X,\,y\,\oplus\,Y}(\,z\,\oplus\,Z\,)\;\Big)
 $$
 Calculating the traces of the following expressions appearing in
 $R_{x\oplus X,\,y\oplus Y}(z\oplus Z)$ over $x\oplus X$
 \begin{eqnarray*}
  \Big(\qquad\;F_+X\;\wedge\;F_+Y\;\qquad\Big)\,(z\oplus Z)
  &=&
  0\,\oplus\,\Big(\;g(\,F_+Z,\,X\,)\,F_+Y\;-\;g(\,F_+Y,\,Z\,)\,F_+X\;\Big)
  \\
  \Big(x\,(F_+^2+[F_+,F_-])\,Y\,\wedge\,\1\Big)\,(z\oplus Z)
  &=&
  \Big(\,g(\,(F_+^2+[F_+,F_-])Y,\,Z\,)\,x\,\Big)\,\oplus\,
  \Big(\quad\;\;\ldots\;\;\quad\Big)
  \\
  \Big(y\,(F_+^2+[F_+,F_-])\,X\,\wedge\,\1\Big)\,(z\oplus Z)
  &=&
  \Big(\qquad\ldots\qquad\Big)\,\oplus\,\Big(\,-yz\,(F_+^2+[F_+,F_-])X\,\Big)
 \end{eqnarray*}
 we obtain $g(\,F_+^2Y\,-\,(\tr\,F_+)\,F_+Y,\,Z\,)$ and $g(\,F_+^2Y\,+\,
 [F_+,F_-]Y,\,Z\,)$ as well as $-yz\,(\tr\,F_+^2)$ so that the Ricci
 endomorphism $g(\,\mathrm{Ric}(\,y\,\oplus\,Y\,),\,z\,\oplus\,Z\,)\,:=\,
 \mathrm{Ric}(\,y\,\oplus\,Y,\,z\,\oplus\,Z\,)$ reads:
 \begin{equation}\label{ric}
  \mathrm{Ric}(\;y\,\oplus\,Y\;)
  \;\;=\;\;
  -\;\Big(\,\tr(\;F_+^2\;)\,y\,\oplus\,
   (\;[\;F_+,\,F_-\;]\;+\;(\tr\,F_+)\,F_+\;)\,Y\;\Big)
 \end{equation}
 In this way we obtain the following upper and lower bound on the stabilizer
 of $R$ in $\so\,V$
 \begin{equation}\label{ulb}
  \stab\;\mathrm{Ric}
  \;\;\supset\;\;
  \stab\;R
  \;\;\supset\;\;
  \stab\;F_+\;\cap\;\stab\;[\,F_+,\,F_-\,]\;\cap\;\so\,V_\circ
 \end{equation}
 because the explicit formula for $R$ tells us that every $\X\,\in\,
 \so\,V_\circ$ stabilizing $F_+$ and $[\,F_+,\,F_-\,]$ stabilizes $R$.
 In order to proceed we need to be somewhat more specific about the
 special form of the endomorphism $F$ in the family of examples found
 by Meusers:

 \begin{Definition}[Meusers' Family of Examples]
 \hfill\label{meu}\break
  An endomorphism $F:\,V_\circ\longrightarrow V_\circ$ on a euclidian
  vector space $V_\circ$ of dimension $m-1\,\geq\,3$ is called special
  provided its diagonalizable symmetric part $F_+$ has only two different
  eigenvalues of multiplicities $1$ and $m-2$ respectively and every
  eigenvector $e\,\neq\,0$ in the $1$--dimensional eigenspace is
  cyclic for the skew symmetric part $F_-$ of $F$ in the sense:
  $$
   V_\circ
   \;\;=\;\;
   \span\;\{\quad e,\;\;F_-\,e,\;\;F_-^2\,e,\;\;F_-^3\,e,\;\;
   F_-^4\,e,\;\ldots\qquad\}
  $$ 
 \end{Definition}

 \noindent
 Using the cyclicity of the eigenvector $e\,\neq\,0$ of the symmetric part
 $F_+$ of a special endomor\-phism $F$ under its skew symmetric part $F_-$
 we may construct a complete flag on $V_\circ$ via:
 $$
  \{\,0\,\}
  \;\subsetneq\;
  \span\;\{\;e\;\}
  \;\subsetneq\;
  \span\;\{\;e,\;F_-e\;\}
  \;\subsetneq\;\ldots\;\subsetneq\;
  \span\;\{\;e,\;F_-e,\;\ldots,\;F_-^{m-3}e\;\}
  \;\subsetneq\;
  V_\circ
 $$
 Up to the choice of signs there exists a unique orthonormal basis $e_2,
 \,\ldots,\,e_m$ of $V_\circ$ adapted to this flag, in this basis the
 matrix of the special endomorphism $F$ is tridiagonal of the form
 \begin{equation}\label{tri}
  F
  \;\;\widehat=\;\;
  \pmatrix{
    f_1 & -f_3 &        &        &       \cr
   +f_3 &  f_2 & -f_4   &        &       \cr
        & +f_4 & \ddots & \ddots &       \cr
        &      & \ddots &   f_2  &  -f_m \cr
        &      &        &  +f_m  &   f_2
  }
 \end{equation}
 where the parameters $f_1,\,\ldots,\,f_m\,\in\,\R$ are arbitrary except
 for $f_1\,\neq\,f_2$ and $f_3,\,\ldots,\,f_m\,\neq\,0$ to ensure the
 cyclicity of the basis vector $e_2$. The special endomorphisms form in
 this way an $m$--parameter family of orbits in $\End\,V_\circ$ under
 the action of the orthogonal group $\mathbf{O}(\,V_\circ\,)$. Extending
 this orthonormal basis to the orthonormal basis $\1,\,e_2,\,\ldots,
 \,e_m$ of $V$ we obtain the following explicit matrix for the Ricci
 endomorphism calculated in equation (\ref{ric}):
 $$
  \mathrm{Ric}
  \;\;\widehat=\;\;
  \pmatrix{
   -\,\tr\,F_+^2 & & & & &             \cr
   & -(\;\tr\,F_+\;)f_1 & (f_1-f_2)f_3 & & & \cr
   & (f_1-f_2)f_3 & -(\;\tr\,F_+\;)f_2 & & & \cr
   & & & -(\,\tr\,F_+\,)f_2 & &              \cr
   & & & & \ddots &                          \cr
   & & & & & -(\,\tr\,F_+\,)f_2
  }
 $$
 This matrix certainly has $4$ different eigenvalues of multiplicities
 $m-3$ and $1,\,1,\,1$ for a generic special endomorphism $F$, say
 for $f_1\,=\,1$ and $f_2\,=\,0$ these eigenvalues are $0$ and $-1,
 \,-\frac12\pm\frac12\sqrt{4f_3^2+1})$ respectively. For such a generic
 special endomorphism $F$
 $$
  \stab\;\mathrm{Ric}
  \;\;=\;\;
  \stab\;F_+\;\cap\;\stab\;[\,F_+,\,F_-\,]\;\cap\;\so\,V_\circ
  \;\;=\;\;
  \so\,\{\;\1,\;e_2,\;e_3\;\}^\perp
 $$
 and equation (\ref{ulb}) tells us that this agrees with the stabilizer
 of $R$:
 $$
  \h_0\;\;:=\;\;\stab\;R\;\;=\;\;\so\,\{\;e_4,\;e_5,\;\ldots,\;e_m\;\}
 $$
 Because every $\X\,\in\,\h_0$ commutes with $F_+$, the comodule directional
 derivatives are given by
 $$
  \frac{\partial\X}{\partial\1}
  \;\;=\;\;
  [\;\X,\;F_-\;]
  \qquad\qquad
  \frac{\partial\X}{\partial X}
  \;\;=\;\;
  [\;\X,\;F_+X\,\wedge\,\1\;]\;-\;
  F_+(\,\X\,X\,)\,\wedge\,\1
  \;\;=\;\;
  0
 $$
 for all $X\,\in\,V_\circ$, in particular for $\X\,:=\,e_\mu\wedge e_\nu\,
 \in\,\h_0$ with $4\leq\mu<\nu\leq m$ we obtain this way
 $$
  \frac\partial{\partial\mathbf{1}}{(\,e_\mu\,\wedge\,e_\nu\,)}
  \;\;=\;\;
  -\;\Big(\;F_-e_\mu\,\wedge\,e_\nu\;+\;e_\mu\,\wedge\,F_-e_\nu\;\Big)
  \;\;\equiv\;\;
  f_\mu\,e_{\mu-1}\,\wedge\,e_\nu
 $$
 modulo $\so\,\{e_\mu,\,\ldots,\,e_m\}$. Consequently the subalgebras
 in the stabilizer filtration (\ref{df}) read
 \begin{equation}\label{mres}
  \h_r
  \;\;=\;\;
  \so\,\{\;e_{r+4},\,\ldots,\,e_m\;\}
 \end{equation}
 for all $r\,=\,0,\ldots,m-4$ due to our cyclicity assumption $f_3,\,
 \ldots,\,f_m\,\neq\,0$. In consequence the affine homogeneous spaces
 of dimension $m$ associated to generic special endomorphisms $F:\,V_\circ
 \longrightarrow V_\circ$ on euclidian vector spaces $V_\circ$ of dimension
 $m-1$ have Singer invariant $m-4$, because $\h_{m-5}\,\neq\,\{0\}$, but
 $\h_{m-4}\,=\,\h_\infty\,=\,\{0\}$.

 Before closing this section we want to calculate the Spencer cohomologies
 in this family of examples of affine spaces with large Singer invariant.
 With all directional derivatives $\frac{\partial\X}{\partial X}\,=\,0$
 in directions $X\,\in\,V_\circ$ vanishing we may use the isomorphism
 (\ref{sc}) to reduce this calculation to the calculation of the Spencer
 cohomology of the comodule $\h^\bullet$ considered as a comodule over
 $\S\,V/V_\circ\,=\,\S\,\R$. By definition the comodule $\h^\bullet$ is
 the direct sum of the successive quotients $\h_{r-1}/\h_r$, which are
 spanned for $r\,>\,0$ by the bivectors $e_{r+3}\wedge e_\nu\,+\,\h_r$
 with $r+3\,<\,\nu\,\leq\,m$. The only non--trivial directional derivatives
 are injective for $r\,>\,0$
 $$
  \frac\partial{\partial\1}:
  \quad\h_r/\h_{r+1}\;\longrightarrow\;\h_{r-1}/\h_r,
  \qquad e_{r+4}\wedge e_\nu\;+\;\h_{r+1}\;\longmapsto\;
  f_{r+3}\,e_{r+3}\wedge e_\nu\;+\;\h_r
 $$
 with cokernel spanned by $e_{r+3}\,\wedge\,e_{r+4}\,+\,\h_r$. In summary
 we have proved the following lemma:

 \begin{Lemma}[Singer Invariant and Spencer Cohomology]
 \hfill\label{scex}\break
  Consider a special endomorphism $F:\,V_\circ\longrightarrow V_\circ$
  on a euclidian vector space $V_\circ$ of dimension $m-1\,\geq\,3$ in
  the sense of Definition \ref{meu} and let $e_2,\,\ldots,\,e_m$ be the
  essentially unique orthonormal basis of $V_\circ$, in which $F$ takes
  the trilinear form (\ref{tri}). The euclidian Lie algebra $\g\,=\,
  \R\oplus V_\circ\,=\,V$ associated to the special endomorphism $F$
  has Singer invariant
  $$
   \Singer(\;\g\;)\;\;=\;\;m\;-\;4
  $$
  provided the Ricci endomorphism of $\g$ has $4$ different eigenvalues.
  For such a generic special endomorphism $F$ the Spencer cohomology
  $H^{r,\circ}(\,\h\,)$ of the associated comodule $\h^\bullet$ over
  the coalgebra $\S\,V^*$ is a free $\L^\circ V^*_\circ$--module of rank
  $1$ for all $r\,=\,1,\,\ldots,\,m-4$ with isomorphism
  $$
   \L^{\circ-1}V^*_\circ\;\stackrel\cong\longrightarrow
   H^{r,\,\circ}(\;\h\;),\qquad
   \eta\;\longmapsto\;\Big[\;\1^\sharp\,\wedge\,\eta\;\otimes\;
   (\,e_{r+3}\wedge e_{r+4}\,+\,\h_r\,)\;\Big]
  $$
  where $\1^\sharp\,\in\,V^*$ is the linear form $\1^\sharp(x\oplus X)
  \,:=\,x$. In particular the interesting Spencer cohomologies $H^{r,1}(\h)$
  for $r\,=\,1,\,\ldots,\,m-4$ are all one--dimensional.
 \end{Lemma}

 \pfill
 The remarkable conclusion of our calculation of the Spencer cohomology
 of Meusers' examples of Riemannian homogeneous spaces with large Singer
 invariant is that this specific family of examples is in a sense maximal
 in the moduli space $\M(\,\so\,V\,)_\infty$: Every vector formally tangent
 to the moduli space $\M(\,\so V\,)_\infty$ in a point corresponding to a
 Riemannian homogeneous space in Meusers' family is integrable to a real
 deformation by changing one of the parameters $f_5,\,\ldots,\,f_m-1\,\neq\,0$
 of this family.
\end{document}